\documentclass[10pt, a4paper]{article}
\usepackage{mathtext}
\usepackage{mathrsfs}
\usepackage{amssymb}
\usepackage{amsfonts}
\usepackage{tensor}
\usepackage{amscd}
\usepackage{amsmath, amsthm}
\usepackage{mathtools}
\usepackage{hyperref}
\usepackage[nice]{nicefrac}
\usepackage[symbol]{footmisc}
\usepackage{empheq}
\usepackage{relsize}

\renewcommand{\section}[2]{}

\newcommand*{\Gp}{\mathrm{G_{P}^{\vphantom{-1}}}}
\newcommand*{\Up}{\mathrm{U_{P}^{\vphantom{-1}}}}
\newcommand*{\Dp}{\mathrm{D_{P}^{\vphantom{-1}}}}
\newcommand*\widebox[1]{\fbox{\hspace{1em}#1\hspace{1em}}}

\begin{document}
\begin{center}\textbf{ON A CERTAIN TYPE OF DEFORMATION OF UMBRAL ORTHOGONAL POLYNOMIALS}
\end{center}
\begin{center}
           \begin{footnotesize}
	DANIL KROTKOV
	\end{footnotesize}
\end{center}

\begin{footnotesize}
\noindent\textsc{Abstract}. In the present paper we derive complicated families of orthogonal polynomials in one variable from scratch using the known ones as building blocks. We recall the basics of operational formalism and introduce the notations we use throughout this paper. Then we establish through this formalism the known families of orthogonal polynomials in order to see the mechanism that lies at the heart of classical and hypergeometric $q=1$ orthogonal polynomials and their associated families, which arise form this viewpoint in a very natural manner.
\end{footnotesize}\\

\begin{center}
\textbf{Preliminaries}
\end{center}

\textbf{Key Definitions}
\begin{itemize}
\small{\item We call a sequence of polynomials $\mathrm{P}=\{p_n(x)\}_{n\in\mathbb{N}_{0}}$ a family of monic polynomials, if $\mathrm{deg}~p_n(x)=n$ and the leading coefficient of each polynomial is equal to $1$. We assume that all polynomials and formal power series are defined over $\mathbb{C}$.\\

\item \noindent For given family of monic polynomials $\mathrm{P}$ there exists linear operator $\mathrm{G_P}$, which acts on polynomials by the rule
$$
\mathrm{G_P} \cdot x^n = p_n(x)
$$
\noindent This operator is invertible, since every monic polynomial may be expressed as a linear combination of $p_k(x)$ uniquely. Explicit computation of these operators is of crucial importance, since they contain all of the information about the properties of the corresponding families.\\

\item \noindent For given sequence $\{ h_n\}_{n\in\mathbb{N}_{0}}$ define an operator $h_\theta$: $h_\theta \cdot x^n = h_n x^n$. If $h_n=h(n)$ for some polynomial $h(x)$, then the following relation holds: $h_\theta=h(\theta)=h(x\mathrm{D})$, where $\mathrm{D}=\frac{d}{dx}$. The general relations are
$$
h_\theta x = x h_{\theta+1} ~~~~~~~~~~ \mathrm{D} h_{\theta} = h_{\theta+1} \mathrm{D} ~~~~~~~~~~ \mathrm{D}x-x\mathrm{D}=1
$$
$$
\ell(\mathrm{D})x=\ell'(\mathrm{D})+x\ell(\mathrm{D})
$$
\item \noindent It is natural to consider ''multiplication'' and ''derivative'' operators
\begin{align*}
\mathrm{U_P}\coloneqq\mathrm{G_{P}^{\vphantom{-1}}}x\mathrm{G}_{\mathrm{P}}^{-1}~~~~~&\Rightarrow~~~~~\mathrm{U_P}\cdot p_n(x)=p_{n+1}(x)\\
\mathrm{D_P}\coloneqq\mathrm{G_{P}^{\vphantom{-1}}}\mathrm{D}\mathrm{G}_{\mathrm{P}}^{-1}~~~~~&\Rightarrow~~~~~\mathrm{D_P}\cdot p_n(x)=np_{n-1}(x)
\end{align*}
with an obvious relation
$$
\mathrm{D_P}\mathrm{U_P}-\mathrm{U_P}\mathrm{D_P}=1
$$
and
$$
\mathrm{U}_\mathrm{P}^n \cdot 1 = p_n(x)
$$
In addition, for given sequence $\{ a_n\}_{n\in\mathbb{N}_{0}}$ define an operator $a_{\theta_\mathrm{P}}\coloneqq \mathrm{G_{P}^{\vphantom{-1}}}a_{\theta}\mathrm{G}_{\mathrm{P}}^{-1}$. Obviously this operator has the property $a_{\theta_{\mathrm{P}}}\cdot p_n(x)=a_np_n(x)$.\\

\item \noindent For given family $\mathrm{P}$ there exists a set of formal power series $\{f_n(y)\}_{n\in\mathbb{N}_{0}}$, s.t. $f_n(y) \in y^n+y^{n+1}\mathbb{C}[[y]]$, and such that the following relation holds
$$
e^{xy}=\sum_{n=0}^\infty \frac{p_n(x)f_n(y)}{n!}
$$
The series $f_n$ may be written down explicitly by the use of polynomials, which belong to the dual family $\mathrm{P}^*=\{p_n^*(x)\}_{n\in\mathbb{N}_{0}}$ of polynomials defined by the relation $p_n^*(x)\coloneqq\mathrm{G}_{\mathrm{P}}^{-1}\cdot x^n$.\\

\item There is an important class of polynomials of binomial type. The polynomials of a given family in this class satisfy the binomial theorem condition
$$
p_n(\alpha+\beta)=\sum_{k=0}^n \left(\vphantom{\bigg|}\genfrac{}{}{0pt}{0}{n}{k}\right) p_k(\alpha)p_{n-k}(\beta) ~~~~~\Rightarrow~~~~~\exists f \in y+y^2\mathbb{C}[[y]]:~~~~~ e^{xy}=\sum_{n=0}^\infty \frac{p_n(x)f^n(y)}{n!}
$$
So that this family is defined by a single invertible power series. The binomial case is so important that its corresponding $\mathrm{G_P}$ operator deserves its own name, we call it $\mathrm{C}_f$. So that
$$
\mathrm{U_P}=\mathrm{C}_f x \mathrm{C}_f^{-1} = x\frac{1}{f'(\mathrm{D})}~~~~~~~~~~~~~~~~\mathrm{D_P}=\mathrm{C}_f \mathrm{D}\mathrm{C}_f^{-1} = f(\mathrm{D})
$$
and we have the following property of umbral composition
$$
\mathrm{C}_f \mathrm{C}_g = \mathrm{C}_{g\circ f}
$$
Also we use the following notation: $\varphi(y)=f^{inv}(y)$ means that $f(\varphi(y))=y$, $\mathfrak{T}f(y) \coloneqq f(y)/f'(y)$ and we usually write $(f/f')^{inv}=\omega$. Notice that by $\mathfrak{T}^2 f$ we mean $ff'/(f'^2-ff'')=\mathfrak{T}(\mathfrak{T}f)$ etc.

\item \noindent There is a much more important class of general orthogonal polynomials. These polynomials are defined by the two sequences $\{a_n\}_{n\in\mathbb{N}_{0}}$, $\{b_n\}_{n\in\mathbb{N}_{0}}$, $b_i\neq 0$ and the corresponding three-term relation
$$
\mathrm{U_P}=x-a_{\mathrm{U_P}\mathrm{D_P}}^{\vphantom{-1}}-\mathrm{D_P}b_{\mathrm{U_P}\mathrm{D_P}}^{\vphantom{-1}} ~~~\iff~~~ \mathrm{U_{P^{\text{*}}}}=x+a_\theta+\mathrm{D}b_\theta
$$
(see \cite{Gdsl}, \cite{AW}). These polynomials are also determined by a single power series, but in this case it is $f_0(y)$. Notice that we do not require any positivity conditions on $\{b_n\}_{n\in\mathbb{N}_{\geqslant 0}}$.\\

\item \noindent For an arbitrary operator $\mathrm{T}_x$, which acts on $\mathbb{C}[x]$, we define an operator $\mathrm{\overline{T}}_y$, which acts on $\mathbb{C}[[y]]$, such that
$$
\mathrm{T}_x\cdot e^{xy}=\mathrm{\overline{T}}_y \cdot e^{xy}
$$
Such an operator always exists and is fully determined by its images on monomials, since for an arbitrary $\ell(y) \in \mathbb{C}[[y]]$ holds
$$
\mathrm{\overline{T}}_y \cdot \ell(y)= \mathrm{\overline{T}}_y \cdot \ell(y)e^{xy}\big|_{x=0}=\mathrm{\overline{T}}_y \ell(\mathrm{D})\cdot e^{xy}\big|_{x=0}=\ell(\mathrm{D})\mathrm{\overline{T}}_y\cdot e^{xy}\big|_{x=0}=\sum_{k=0}^\infty \frac{y^k}{k!} [\ell(\mathrm{D})\mathrm{T}_x \cdot x^n]\big|_{x=0}
$$
$\mathrm{T}_x$ sends polynomials to polynomials, so the action of $\ell(\mathrm{D})$ is well-defined and thus the action of $\mathrm{\overline{T}}_y$ is well-defined too. Note that the action of $\mathrm{T}_x$ can be reconstructed from the action of $\mathrm{\overline{T}}_y$. We denote $\mathrm{\overline{U}_P}$ by $\mathfrak{d}_{\mathrm{P}}$ and $\mathrm{\overline{D}_P}$ by $\mathfrak{u}_{\mathrm{P}}$. In the case of binomial monic family we have $\mathrm{D_P}=f(\mathrm{D})$, hence $\mathfrak{u}_{\mathrm{P}}=f(y)$, and $\mathrm{U_P}=xf'(\mathrm{D})^{-1}$, hence $\mathfrak{d}_{\mathrm{P}}=\frac{d}{df}$. Finally we define $\overline{\Gp}=\mathfrak{g}_{\mathrm{P}}^{-1}$ and exceptional operator $\overline{\mathrm{C}_f}=\mathfrak{c}_f^{-1}$ for binomial family, associated to $f^n(y)$ Notice that
$$
\mathfrak{c}_f \cdot y^n = f(y)^n
$$
$$
\mathfrak{g}_{\mathrm{P}} \cdot y^n = f_n(y)
$$
For diagonal operators we do not change the notation: $\overline{\theta_\mathrm{P}}=\theta_\mathrm{P}$, $\overline{a_{\theta_\mathrm{P}}}=a_{\theta_\mathrm{P}}$. Also obviously
$$
\mathfrak{d}_{\mathrm{P}}\mathfrak{u}_{\mathrm{P}}-\mathfrak{u}_{\mathrm{P}}\mathfrak{d}_{\mathrm{P}}=1
$$
\item \noindent Given two different raising operators $\mathrm{U_P}$ and $\mathrm{U_Q}$, their affine combination $t\mathrm{U_P}+(1-t)\mathrm{U_Q}$ is also a raising operator. And more generally, given in addition an invertible sequence $\mathcal{F}_\theta$, $\mathcal{F}_0=1$, $\mathcal{F}_{\theta+1}/\mathcal{F}_\theta = F_\theta$, we may consider an operator $t\mathrm{U_P}F_\theta+(1-t)\mathrm{U_Q}$ for formal variable $t$. With a great abuse of language the process of finding such an invertible operator $\mathcal{T}_t$, that $\mathcal{T}_t x \mathcal{T}_t^{-1} =$ $= t\mathrm{U_P}F_\theta+(1-t)\mathrm{U_Q}$ we call a diagonalization. The corresponding polynomials, now not necessary monic may then be expressed as $\mathcal{T}_t \cdot x^n = (t\mathrm{U_P}F_\theta+(1-t)\mathrm{U_Q})^n \cdot 1 $
}\end{itemize}
\noindent\rule{\textwidth}{1.2pt}
\newpage
\noindent [Binomial families overview]. Suppose we have monic family which satisfies the binomial theorem condition
$$
p_n(\alpha+\beta)=\sum_{k=0}^n \binom{n}{k}p_k(\alpha)p_{n-k}(\beta) ~~~~~\Rightarrow~~~~~\exists f \in x+x^2\mathbb{C}[[x]]:~~~~~ e^{xy}=\sum_{n=0}^\infty \frac{p_n(x)f^n(y)}{n!}
$$
Define an operator $\mathfrak{c}_f:~\mathfrak{c}_f \cdot g(x)=g(f(y))$. Define in addition $C_f = (\overline{\mathfrak{c}}_f)^{-1}$ where $C$ stands for composition. Consider now the functional inverse to $f^{inv}=\varphi$. We then have the following identities
$$
\Gp = \mathrm{C}_f ~~~~~~~~~~  \Up = \mathrm{A}_f \coloneqq x\frac{1}{f'(\mathrm{D})} ~~~~~~~~~~ \Dp = \mathrm{D}_f \coloneqq f(\mathrm{D})
$$
$$
\overline{\Gp}=\mathfrak{g}_{\mathrm{P}}^{-1}=\mathfrak{c}_f^{-1}=\mathfrak{c}_\varphi ~~~~~~~~~~ \mathfrak{u}_{\mathrm{P}}=y_f=f(y)~~~~~~~~~~\mathfrak{d}_\mathrm{P}=\frac{d}{df} ~~~~~~~~~~ f_n^*(y)=\varphi^{n}(y)
$$
Now notice that the classical Lagrange inversion formula implies the following identity
$$
p_n(x)=x\left(\frac{\mathrm{D}}{f(\mathrm{D})}\right)^n \cdot ~x^{n-1}=f'(\mathrm{D})\left(\frac{\mathrm{D}}{f(\mathrm{D})}\right)^{n+1} \cdot ~x^n
$$
The form of this expression allows us to extend the definition of $p_n$ to an arbitrary complex index and moreover to an arbitrary parameter $t$ in the following manner
$$
p_s^{t}(\alpha)=\alpha \left(\frac{\mathcal{D}f'(t)}{f(\mathcal{D}+t)-f(t)}\right)^{s}\cdot~\alpha^{s-1}
$$
Notice that one could think of $t$ as of formal point in complex plane, meaning that if the formal series $f$ defines an entire function with simple zeros we could substitute $t=z$ for arbitrary $z \in \mathbb{C}$. The resulting series is of the form $\alpha^s+\alpha^{s-1}\mathbb{C}[[t]][[\alpha^{-1}]]$ and satisfies $\mathcal{A}_f \cdot p_s^{t}(\alpha)f'(t)^{-s}e^{\alpha t} =$ $= p_{s+1}^{t}(\alpha)f'(t)^{-s-1}e^{\alpha t}$ since there is an equality of operators $g(\mathcal{D})e^{\alpha t}=e^{\alpha t}g(\mathcal{D}+t)$. This property allows us to define the complex powers of $\mathcal{A}_f$ as follows
$$
\mathcal{A}_f^s=(\alpha)p_s^{\mathcal{D}}f'(\mathcal{D})^{-s}
$$
Notice that substitution $s=n$ recovers the formula for $(\mathrm{U}_{\mathrm{P}})^n$. Indeed $\mathcal{D}_f\mathcal{A}_f^s=s\mathcal{A}_f^{s-1}+\mathcal{A}_f^s\mathcal{D}_f$ and $\mathcal{A}_f^{s+h}=\mathcal{A}_f^s\mathcal{A}_f^h$. This operator has also the following representation, which also degenerates to the usual one when $s=n$.
$$
\mathcal{A}_f^s=(\alpha)\frac{p_{s+\mathcal{A}_f\mathcal{D}_f}}{p_{\mathcal{A}_f\mathcal{D}_f}}
$$
The latter equality of operators leads to the limit formula, which resembles the property of Harmonic numbers and factorials. Consider the series $(f/f')^{inv}=\omega$. Then the following holds
\begin{align*}
\lim_{s\to\infty} \frac{\dot p_{s}^{\vphantom{\omega}}(s\alpha)}{p_{s}(s\alpha)}-\ln s=\ln\alpha-\ln f'(\omega(\alpha^{-1}))~~~~~~~~\lim_{s\to\infty} \frac{p_{s+\text{\scriptsize{H}}}^{\vphantom{\omega}}(t+s\alpha)}{p_{s}(s\alpha)}s^{-\text{\scriptsize{H}}}=f(\omega(\alpha^{-1}))^{-\text{\scriptsize{H}}}e^{t\omega(\alpha^{-1})}
\end{align*}
It also leads to the following logarithmic version of Lagrange inversion formula
$$
\frac{1}{s}\ln p_s(\alpha)=\exp\left(\frac{\partial}{\partial\alpha}\left(\frac{d}{d\omega}-s\mathrm{L}\right)\vphantom{\bigg|}\right)\cdot~\frac{\omega(x)}{x}\bigg|_{x=0} \cdot ~\ln \alpha
$$
which can be used to calculate the formal asymptotics of $\ln p_s(s\alpha^{-1})$ for large $s$ (so that it could be seen as a suitable replacement of the Euler-Maclaurin summation formula, which is commonly used for this purpose in particular case $f(x)=e^x-1$ and associated ones):
\begin{subequations}
\begin{empheq}[box=\widebox]{align*}
&~\\
&\ln p_s(s\alpha^{-1})\sim s\ln (s\alpha^{-1})-s\alpha^{-1}\int_{0}^\alpha\ln f'(\omega(t))dt+\frac{1}{2}\ln\omega'(\alpha)+\\
&~~~~~+\frac{1}{24s}\left(2\frac{\omega'(\alpha)-1}{\omega'(\alpha)}+4\frac{\alpha^2\omega''(\alpha)^2}{\omega'(\alpha)^3}-2\frac{\alpha\omega''(\alpha)}{\omega'(\alpha)^2}-3\frac{\alpha^2 \omega^{(3)}(\alpha)}{\omega'(\alpha)^2}\right)-\frac{1}{48s^2}\frac{\alpha^{3}}{\omega'(\alpha)}\left(\frac{\alpha}{\omega'(\alpha)}\right)^{(4)}+...\\
\end{empheq}
\end{subequations}
~\\
\noindent\rule{\textwidth}{1pt}
~\\
\newpage
\noindent [Orthogonal families overview]. We now list the basics of algebraic aspects of continued fractions and orthogonal polynomials. Suppose we have the following continued fraction of formal power series $x F_0(x)\in x+x^2\mathbb{C}[[x]]$, ($b_i \neq 0$).
$$
xF_0(x)=\cfrac{x}{1-a_0x-\cfrac{1\cdot b_1 x^2}{1-a_1x-\cfrac{2\cdot b_2 x^2}{1-a_2x-\cfrac{3\cdot b_3 x^2}{1-a_3x-...}}}}
$$
In order to find the $n$-th convergent we consider the following product of matrices
$$
\left( \begin{array}{cc}
0 & x \\
-b_1 x & 1-a_0 x
\end{array} \right)
\left( \begin{array}{cc}
0 & x \\
-2b_2 x & 1-a_1 x
\end{array} \right)\cdot...\cdot
\left( \begin{array}{cc}
0 & x \\
-nb_n x & 1-a_{n-1} x
\end{array} \right)=
\left( \begin{array}{cc}
-nb_n x R_{n-1} & R_n \\
-nb_n x Q_{n-1} & Q_n
\end{array} \right)
$$
We then have the following three-term recurrence relations (we assume that $b_0 \in \mathbb{C}$ is an arbitrary quantity and it does not give any contribution to the recurrence relation).
\begin{align*}
R_{n+1}(x)&=R_n(x)(1-a_n x)-nb_nx^2R_{n-1}(x) ~~~~~~~~~~~R_0=0,~~R_1=x\\
Q_{n+1}(x)&=Q_n(x)(1-a_n x)-nb_nx^2Q_{n-1}(x) ~~~~~~~~~~~Q_0=1,~~Q_1=1-a_0x
\end{align*}
Now introduce the monic family $p_n(x)=x^nQ_n(x^{-1})$ with $\mathrm{U_P}=x-a_{\mathrm{U_P}\mathrm{D_P}}^{\vphantom{-1}}-\mathrm{D_P} b_{\mathrm{U_P}\mathrm{D_P}}^{\vphantom{-1}}$  the corresponding raising operator. Denote $b_1b_2...b_n=\mathcal{B}_n$. Then the equality of determinants corresponds to the series
$$
 \frac{R_{n+1}}{Q_{n+1}}-\frac{R_{n}}{Q_{n}}=\frac{n!\mathcal{B}_nx^{2n+1}}{Q_nQ_{n+1}}~~~~~\Rightarrow~~~~~x^{-1}F_0(x^{-1})=\sum_{n=0}^\infty \frac{n!\mathcal{B}_n}{p_n(x)p_{n+1}(x)}
$$
Consider now the formal power series $f_0(x)=\theta!^{-1}F_0(x)$, which is the inverse Laplace transform of $F_0$. Define general inner product of polynomials by
$$
\langle h_1,h_2 \rangle_{f_0}=f_0(\mathrm{D}) h_1h_2\big|_{x=0}
$$
Then the condition $F_0(x)-(xQ_n)^{-1}R_n=O(x^{2n})$ implies that polynomials $p_n$ are orthogonal with respect to this inner product and the following formula for their norm is valid
$$
\langle p_n,p_n\rangle_{f_0}=n!\mathcal{B}_n
$$
The corresponding family $f_n(y)$ with corresponding derivative operator $\mathfrak{d}_{\mathrm{P}}=\mathfrak{d}-a_{\mathfrak{u}_\mathrm{P}\mathfrak{d}_{\mathrm{P}}}-b_{\mathfrak{u}_\mathrm{P}\mathfrak{d}_{\mathrm{P}}}\mathfrak{u}_\mathrm{P}$ is then given by
$$
e^{xy}=\sum_{n=0}^\infty \frac{p_n(x)f_n(y)}{n!}~~; ~~~~~~~~~~ \langle p_n, e^{xy}\rangle_{f_0}=\mathcal{B}_nf_n(y)~~~~~\Rightarrow~~~~~p_n(\mathfrak{d})\cdot f_0(y)=\mathcal{B}_nf_n(y)
$$
The latter equalities imply the addition theorem with formula for $n!F_n(x)=\theta!f_n(x)$ as consequence
$$
f_0(y+t)=e^{t\mathfrak{d}}\cdot f_0(y) ~\Rightarrow~ f_0(y+t)=\sum_{n=0}^\infty \frac{\mathcal{B}_n}{n!}f_n(y)f_n(t) ~\Rightarrow~ \frac{xF_0(x)-yF_0(y)}{x-y}=\sum_{n=0}^\infty n!\mathcal{B}_nF_n(x)F_n(y)
$$
The three-term recurrence relation for $p_n$ imply Christoffel-Darboux formula
$$
\frac{p_n(y)p_{n+1}(x)-p_{n+1}(y)p_n(x)}{x-y}=\sum_{k=0}^{n}\frac{n!\mathcal{B}_n}{k!\mathcal{B}_k}p_k(x)p_k(y)
$$
and these polynomials, viewed as polynomials in $x$ parametrized with $y$, are orthogonal too. The corresponding $\mathrm{G_{Q}}$-operator is given by
$$
\mathrm{G_{Q}}=\mathrm{G_P}\frac{p_\theta(y)}{\mathcal{B}_\theta}(1-\mathrm{D})^{-1}\frac{\mathcal{B}_\theta}{p_\theta(y)}
$$
And its three-term relation in case $y \neq 0$ is given by
\begin{align*}
\mathrm{U_{Q^*}}&=\mathrm{G_{Q}^{-1}} x \mathrm{G_{Q}} = \frac{p_\theta(y)}{\mathcal{B}_\theta} (1-\mathrm{D}) \frac{\mathcal{B}_\theta}{p_\theta(y)} [x+a_\theta + \mathrm{D}b_\theta] \frac{p_\theta(y)}{\mathcal{B}_\theta}(1-\mathrm{D})^{-1}\frac{\mathcal{B}_\theta}{p_\theta(y)}=\\
&= \frac{p_\theta(y)}{\mathcal{B}_\theta} (1-\mathrm{D}) \left[x\frac{b_{\theta+1}p_\theta(y)}{p_{\theta+1}(y)}+a_\theta+\frac{p_{\theta+1}(y)}{p_{\theta}(y)}\mathrm{D}\right](1-\mathrm{D})^{-1}\frac{\mathcal{B}_\theta}{p_\theta(y)}=\\
&= \frac{p_\theta(y)}{\mathcal{B}_\theta}\left[x\frac{b_{\theta+1}p_\theta(y)}{p_{\theta+1}(y)}+a_{\theta+1}-\frac{p_{\theta+1}(y)}{p_{\theta}(y)}+\frac{p_{\theta+2}(y)}{p_{\theta+1}(y)}(1+\mathrm{D})\right]\frac{\mathcal{B}_\theta}{p_\theta(y)}=\\
&= x+a_{\theta+1}-\frac{p_{\theta+1}(y)}{p_{\theta}(y)}+\frac{p_{\theta+2}(y)}{p_{\theta+1}(y)}+\frac{p_{\theta+2}(y)p_\theta(y)b_{\theta+1}}{p_{\theta+1}(y)^2}\mathrm{D}
\end{align*}
Notice that that the Christoffel-Darboux deformation is the deformation of the given $\mathrm{G_P}$ from the right (and it is easier to think of such transformations as of "geometric" ones). In this paper we are interested in the left transformations that are more algebraic in their nature. Obviously the deformations from the right \textit{are} the deformations from the left after conjugation, but it is always hard to understand what the resulting operator actually is. For example how it acts on a the trivial family $\{x^n\}_{n\in\mathbb{N}_{\geqslant{0}}}$. So we distinguish between these types of deformations only informally. To give a glimpse on deformationas that we call \textit{algebraic}, recall that the numerator polynomials may be recovered as
$$
y^n R_n(y^{-1})=f_0(\mathrm{D})\frac{p_n(x)-p_n(y)}{x-y}\bigg|_{x=0}
$$
They satisfy the same recurrence relation but with another initial condition, and we define $p_n(x,1)\coloneqq$ $=x^{n+1} R_{n+1}(x^{-1})$. These polynomials are monic and their recurrence relation is the following one
$$
\mathrm{U_{P^*}(1)}=x+a_{\theta+1}+\frac{(2+\theta)b_{\theta+2}}{1+\theta}\mathrm{D}
$$
Such polynomials are called associated or sometimes corecursive (for further information see \cite{Lewanowicz1},\cite{Lewanowicz2},\cite{ILVW},\cite{Letessier},\cite{W}) and we may define $p_n(x,c)$ for possibly complex values of $c$ (but in general only for nonnegative integral) to be monic orthogonal polynomials with the corresponding dual raising operator
$$
\mathrm{U_{P^*}}(c)=x+a_{\theta+c}+\frac{(1+c+\theta)b_{\theta+c+1}}{1+\theta}\mathrm{D}
$$
The natural question is how to determine its $\mathrm{G_P}$-operator explicitly. Notice that the formula
$$
p_n(x,1)=f_0(\mathrm{D}_y)\frac{p_{n+1}(x)-p_{n+1}(y)}{x-y}\bigg|_{y=0}
$$
implies
$$
\mathrm{G_P}(1)\cdot x^n = \delta_y f_0(\mathrm{D}_y)(1-y\mathrm{L}_x)^{-1}\mathrm{L}_x \cdot p_{n+1}(x)=F_0(\mathrm{L})\mathrm{L}\cdot p_{n+1}(x) ~~\Rightarrow~~ \mathrm{G_P}(1) = F_0(\mathrm{L})\mathrm{L}\mathrm{G_P}x
$$
where $\mathrm{L}$ is the 0-derivative $\mathrm{L}=(1+\theta)^{-1}\mathrm{D}=\theta! \mathrm{D} \theta!^{-1}$ and $\delta=1-x\mathrm{L}$ an operator that takes the value of a polynomial or a power series at zero. Notice also that although we have written $\mathrm{G_P}$ explicitly, we do not know the actual inverse of this operator, since both $\mathrm{L}$ and $x$ are not invertible (altough $\mathrm{L}x \equiv 1$). Now it is easy to see from the construction, that if we define the Laplace transform $F_n(y)\coloneqq n!^{-1}\theta!\cdot f_n(y) \in y^n+y^{n+1}\mathbb{C}[[y]]$ (so that $p_n(\mathrm{L})\cdot F_0(y)=n!\mathcal{B}_nF_n(y)$) then
$$
\frac{F_n(x)}{x F_{n-1}(x)}=\cfrac{1}{1-a_n x-\cfrac{(n+1)\cdot b_{n+1} x^2}{1-a_{n+1}x-\cfrac{(n+2)\cdot b_{n+2} x^2}{1-a_{n+2}x-\cfrac{(n+3)\cdot b_{n+3} x^2}{1-a_{n+3}x-...}}}}
$$
where we naturally set $F_{-1}(y)\coloneqq y^{-1}$ for any monic family (see [Appendix A]). Hence for the associated moment generating function we expect the identity
$$
f_0(y, c)=\theta!^{-1}\cdot \frac{F_c(y)}{y F_{c-1}(y)}
$$
In what follows we establish these deformations of known orthogonal families by the use of operational methods, providing the corresponding $\mathrm{G_P}$-operators explicitly.
\newpage
\noindent [Base case. Sheffer orthogonal family]. Recall that Sheffer families of monic polynomials are defined as families with $\mathrm{G_P}=\ell(\mathrm{D})\mathrm{C}_f$ for series $\ell(y) \in 1+y\mathbb{C}[[y]]$ and $f(y)\in y+y^2\mathbb{C}[[y]]$. There is an easy way to classify all Sheffer orthogonal polynomials. Notice that there exist formal power series $k_1(y), k_2(y)$, such that
$$
\mathrm{G_P^{-1}} x \mathrm{G_P} = x + k_1(\mathrm{D}) + \theta k_2(\mathrm{D})
$$
(since for the functional inverse of $f$, call it $\varphi$ we have $\mathrm{C}_{f}^{-1} \ell(\mathrm{D})^{-1} x \ell(\mathrm{D})\mathrm{C}_{f}=x\varphi'(\mathrm{D})^{-1} -(\ell'/\ell)(\varphi(\mathrm{D}))$ and $x\mathrm{D}=\theta$). Thus we have a restriction on $k_1, k_2$ to be polynomials of degree at most 1. Hence the Sheffer orthogonal polynomials have the following dual raising operator:
$$
\mathrm{U_{P^{*}}}=x+(a+b\theta)+(p+q\theta)\mathrm{D}
$$
In order to simplify all of the computations from now on we assume that the dual raising operator is actually of the form
$$
\mathrm{U_{P^*}}=x+a(1+\lambda \theta)+b(2+\lambda \theta)\mathrm{D}=(1+\lambda a\mathrm{D}+ \lambda b\mathrm{D}^2)^{1/\lambda} x (1+\lambda a\mathrm{D}+\lambda b\mathrm{D}^2)^{1-1/\lambda}
$$
(notice that since we deal with formal algebraic polynomials we do not require $b\neq 0$. Also we work modulo constant term because there always is an elementary deformation $\mathrm{G_Q}=e^{c\mathrm{D}}\mathrm{G_P}$)
That is to say we fix $f(y)$: $f'(y)=1+\lambda af(y)+\lambda bf^2(y)$ and $\ell(y)=f'(y)^{-1/\lambda}$. This is done in such a way so that we may set $\lambda=0$ to obtain f.e. Hermite polynomials. So we fix $\mathrm{G_P}=f'(\mathrm{D})^{-1/\lambda}\mathrm{C}_{f}$. The generating function of the resulting polynomials may thus be obtained as follows
$$
\mathrm{G_P}\cdot e^{xy}=f'(\mathrm{D})^{-1/\lambda}\mathrm{C}_{f} \cdot e^{xy} = f'(\mathrm{D})^{-1/\lambda} \cdot e^{x\varphi(y)} = \varphi'(y)^{1/\lambda} e^{x\varphi(y)}
$$
\noindent\rule{\textwidth}{1.2pt}
~\\
\noindent [First step of deformation. Ultraspherical]. We now raise the following question. Given an invertible operator $\mathcal{F}_\theta$ (invertibility of this operator is equivalent to the condition $\mathcal{F}_i\neq 0$) with $\mathcal{F}_{\theta+1}/\mathcal{F}_\theta = F_\theta$, is there a way to understand the new family $\mathrm{Q}$ with the following dual raising operator?
$$
\mathrm{U_{P^*}}=x+a_\theta+\mathrm{D}b_{\theta} ~~~~\Rightarrow ~~~~\mathrm{U_{Q^*}}=x+a_\theta F_\theta + b_{\theta+1}F_{\theta+1}F_{\theta}\mathrm{D}
$$
This question is difficult to answer in general. What we may do is to seek for the well-understood operator $\mathrm{G}_\mu$ such that
$$
\mathrm{G_Q}=\mathcal{F}_\theta^{-1} \mathrm{G}_\mu \mathrm{G_P} \mathcal{F}_\theta
$$
Then we obtain tautologically the following equation
$$
 \mathrm{G}_\mu^{-1} xF_\theta  \mathrm{G}_\mu = xF_{\theta_\mathrm{P}}
$$
So that we may write $\mathcal{F}_n \mathrm{G}_\mu^{-1}\cdot x^n=\mathcal{F}_n \mu_n^*(x) \coloneqq (xF_{\theta_\mathrm{P}})^n \cdot 1$. But of course it would be nice to obtain the closed form for $\mathrm{G}_\mu^{-1}$ or $\mathrm{G}_\mu$. And in the case of Sheffer family and a specific choice of $\mathcal{F}_\theta$ we may actually write an explicit answer.\\

So we consider the series $f$: $f'(y)=1+\lambda af(y)+\lambda b f^2(y)$, $\varphi=f^{inv} \Rightarrow (\varphi')^{-1}=1+\lambda ay+ \lambda by^2$. Set $\mathrm{G_{P}}=f'(\mathrm{D})^{-1/\lambda}\mathrm{C}_f$, i.e. $\mathrm{G_{P}^{-1}}=(1+\lambda a\mathrm{D}+ \lambda b\mathrm{D}^2)^{1/\lambda}\mathrm{C}_\varphi$. And now consider the following sequence: $\mathcal{F}_\theta = \lambda^{-\theta}\Gamma(\nicefrac{1}{\lambda})\Gamma(\theta+\nicefrac{1}{\lambda})^{-1}~\Rightarrow~ \mathcal{F}_{\theta+1}/\mathcal{F}_\theta=F_\theta=(1+\lambda \theta)^{-1}$. What we have to do is to find the corresponding $\mathrm{G}_\mu$, i.e. solve the equation
$$
 \mathrm{G}_\mu^{-1} x\frac{1}{1+\lambda\theta}  \mathrm{G}_\mu = x f'(\mathrm{D})^{-1/\lambda}\mathrm{C}_f \frac{1}{1+\lambda\theta} \mathrm{C}_f^{-1}f'(\mathrm{D})^{1/\lambda}
$$
Now $\theta=x\mathrm{D}$, hence $\mathrm{C}_f \theta \mathrm{C}_f^{-1} = xf'(\mathrm{D})^{-1}f(\mathrm{D}) = x \mathfrak{T}f(\mathrm{D})$ (recall that we denote $\mathfrak{T}f=f/f'$). Thus the following holds
\begin{align*}
x f'(\mathrm{D})^{-1/\lambda}\frac{1}{1+\lambda x \mathfrak{T}f(\mathrm{D})}f'(\mathrm{D})^{1/\lambda}=x\left(1+\lambda \left(x-\frac{1}{\lambda}\frac{f''}{f'}(\mathrm{D})\right)\mathfrak{T}f(\mathrm{D})\right)^{-1}=\\
x \frac{1}{(\mathfrak{T}f)'(\mathrm{D})+\lambda x \mathfrak{T}f(\mathrm{D})}=x\frac{1}{(\mathfrak{T}f)'(\mathrm{D})}\frac{1}{1+\lambda x (\mathfrak{T}^2 f)(\mathrm{D})}=\mathrm{C}_{\mathfrak{T}f} x \frac{1}{1+\lambda \theta}\mathrm{C}_{\mathfrak{T}f}^{-1}
\end{align*}
And this trick works actually for any Sheffer sequence. So surprisingly we have obtained the following result:
\newpage
Consider the new family $\mathrm{Q}$ with $\mathrm{G_{Q}}=\mathcal{F}_\theta^{-1}\mathrm{C}_{\mathfrak{T}f}^{-1}\mathrm{G_{P}}\mathcal{F}_\theta=\mathcal{F}_\theta^{-1}\mathrm{C}_{\mathfrak{T}f}^{-1}f'(\mathrm{D})^{-1/\lambda}\mathrm{C}_f\mathcal{F}_\theta$. Then its corresponding dual raising operator has the following explicit form
\begin{align*}
\mathrm{U_{Q^*}}&=\mathrm{G^{-1}_{Q}}x\mathrm{G_{Q}}=\mathcal{F}_\theta^{-1}\mathrm{G^{-1}_{P}}\mathrm{C}_{\mathfrak{T}f}\mathcal{F}_\theta [x] \mathcal{F}_\theta^{-1}\mathrm{C}_{\mathfrak{T}f}^{-1}\mathrm{G_{P}}\mathcal{F}_\theta=\\
&=\mathcal{F}_\theta^{-1}\mathrm{G^{-1}_{P}}\mathrm{C}_{\mathfrak{T}f} \left[x\frac{\mathcal{F}_{\theta+1}}{\mathcal{F}_\theta}\right]\mathrm{C}_{\mathfrak{T}f}^{-1}\mathrm{G_{P}}\mathcal{F}_\theta=\mathcal{F}_\theta^{-1}\mathrm{G^{-1}_{P}} \left[x f'(\mathrm{D})^{-1/\lambda}\frac{1}{1+\lambda x \mathfrak{T}f(\mathrm{D})}f'(\mathrm{D})^{1/\lambda}\right]\mathrm{G_{P}}\mathcal{F}_\theta=\\
&=\mathcal{F}_\theta^{-1}\mathrm{G^{-1}_{P}} \left[x f'(\mathrm{D})^{-1/\lambda}\mathrm{C}_f\frac{1}{1+\lambda \theta}\mathrm{C}_f^{-1} f'(\mathrm{D})^{1/\lambda}\right]\mathrm{G_{P}}\mathcal{F}_\theta=\mathcal{F}_\theta^{-1}\mathrm{G^{-1}_{P}} x \mathrm{G_{P}}\frac{1}{1+\lambda \theta}\mathrm{G^{-1}_{P}}\mathrm{G_{P}}\mathcal{F}_\theta=\\
&=\mathcal{F}_\theta^{-1}\mathrm{U_{P^*}}\frac{1}{1+\lambda \theta}\mathcal{F}_\theta=\mathcal{F}_\theta^{-1}\mathrm{U_{P^*}}\mathcal{F}_{\theta+1}=x+a+b\frac{2+\lambda\theta}{(1+\lambda\theta)(1+\lambda+\lambda\theta)}\mathrm{D}
\end{align*}
Notice that in addition we may actually compute the corresponding dual derivative operator
\begin{align*}
\mathrm{D_{Q^*}}&=\mathrm{G^{-1}_{Q}}\mathrm{D}\mathrm{G_{Q}}=\mathcal{F}_\theta^{-1}\mathrm{G^{-1}_{P}}\mathrm{C}_{\mathfrak{T}f}\mathcal{F}_\theta [\mathrm{D}]\mathcal{F}_\theta^{-1}\mathrm{C}_{\mathfrak{T}f}^{-1}\mathrm{G_{P}}\mathcal{F}_\theta=\mathcal{F}_\theta^{-1}\mathrm{G^{-1}_{P}}\mathrm{C}_{\mathfrak{T}f}\left[\frac{\mathcal{F}_\theta}{\mathcal{F}_{\theta+1}}\mathrm{D}\right]\mathrm{C}_{\mathfrak{T}f}^{-1}\mathrm{G_{P}}\mathcal{F}_\theta=\\
&=\mathcal{F}_\theta^{-1}\mathrm{G^{-1}_{P}}\mathrm{C}_{\mathfrak{T}f}[(1+\lambda\theta)\mathrm{D}]\mathrm{C}_{\mathfrak{T}f}^{-1}\mathrm{G_{P}}\mathcal{F}_\theta=\mathcal{F}_\theta^{-1}\mathrm{G^{-1}_{P}}\left[\left(1+\lambda x \frac{\mathfrak{T}f}{(\mathfrak{T}f)'}(\mathrm{D})\right)\mathfrak{T}f(\mathrm{D})\right]\mathrm{G_{P}}\mathcal{F}_\theta=\\
&=\mathcal{F}_\theta^{-1}\mathrm{G^{-1}_{P}}((\mathfrak{T}f)'(\mathrm{D})+\lambda x \mathfrak{T}f(\mathrm{D}))(\mathfrak{T}^2f)(\mathrm{D})\mathrm{G_{P}}\mathcal{F}_\theta=\\
&=\mathcal{F}_\theta^{-1}\mathrm{G^{-1}_{P}}\left[\left(1+\lambda \left(x-\frac{1}{\lambda}\frac{f''}{f'}(\mathrm{D})\right)\mathfrak{T}f(\mathrm{D})\right)(\mathfrak{T}^2f)(\mathrm{D})\right]\mathrm{G_{P}}\mathcal{F}_\theta=\\
&=\mathcal{F}_\theta^{-1}\mathrm{G^{-1}_{P}}[f'(\mathrm{D})^{-1/\lambda}\mathrm{C}_f(1+\lambda\theta)\mathrm{C}^{-1}_f f'(\mathrm{D})^{1/\lambda}(\mathfrak{T}^2f)(\mathrm{D})]\mathrm{G_{P}}\mathcal{F}_\theta=\\
&=\mathcal{F}_\theta^{-1}(1+\lambda\theta)(\mathfrak{T}^2f)(\varphi(\mathrm{D}))\mathcal{F}_\theta=\\
&=\mathcal{F}_{\theta+1}^{-1} \frac{\mathrm{D}}{1-\lambda b \mathrm{D}^2}\mathcal{F}_\theta
\end{align*}
Where in the last line we used the simple observation that if $f'(y)=1+\lambda af(y) +\lambda bf^2(y)$, then $(\mathfrak{T}f)'=(1-\lambda bf^2)/f'$, hence $(\mathfrak{T}^2f)(\varphi)=(\mathfrak{T}f)(\varphi)/(\mathfrak{T}f)'(\varphi)=y\varphi'(y)/((1-\lambda by^2)\varphi'(y))$. The resulting orthogonal polynomials are the classical ultraspherical polynomials with moment generating function
$$
f^{\mathrm{Q}}_0(y)=\overline{\mathrm{G_Q^{-1}}}\cdot 1=\mathfrak{g}_{\mathrm{Q}}\cdot 1=e^{ay} [f^{\mathrm{Q}}_0 |_{a=0}](y)
$$
$\mathrm{G_{Q}}=\mathcal{F}_\theta^{-1}\mathrm{C}_{\mathfrak{T}f}^{-1}\mathrm{G_{P}}\mathcal{F}_\theta$ ~~$\Rightarrow$~~ $[y^n] [f^{\mathrm{Q}}_0 |_{a=0}]=\mathcal{F}_n [y^n] f'(\omega(y))^{1/\lambda}$ ~with~ $f|_{a=0}: f'=1+\lambda bf^2$. Hence
$$
\boxed{
~~f^{\mathrm{Q}}_0(y)=e^{ay}\sum_{n=0}^\infty \frac{(\lambda^{-1})!}{(\lambda^{-1}+n)!} \frac{(b\lambda^{-1})^n y^{2n}}{n!}~~
}
$$
Now $\mathrm{G_Q}$ is given explicitly, hence we may now recover the generating function as $\mathrm{G_{Q}}\mathcal{F}_\theta^{-1} \cdot e^{xy}=$\\ $=\mathcal{F}_\theta^{-1}\mathrm{C}_{\mathfrak{T}f}^{-1}\mathrm{G_{P}} \cdot e^{xy}=\mathfrak{g}_{\mathrm{P}}^{-1}\mathfrak{c}_{\mathfrak{T}f}\mathcal{F}_{\theta_{y}}^{-1} \cdot e^{xy}=(\varphi')^{1/\lambda}\mathfrak{c}_{\varphi}\mathfrak{c}_{\mathfrak{T}f}\mathcal{F}_{\theta_{y}}^{-1} \cdot e^{xy}=(\varphi')^{1/\lambda}\mathfrak{c}_{y\varphi'(y)}\cdot (1-\lambda xy)^{-\frac{1}{\lambda}}$, so that the generating function is equal to
$$
\sum_{n=0}^\infty \binom{-\lambda^{-1}}{n}(-\lambda y)^n q_n(x)=((\varphi'(y))^{-1}-\lambda xy)^{-\frac{1}{\lambda}}=(1+\lambda (a-x)y +\lambda b y^2)^{-\frac{1}{\lambda}}
$$
\noindent\rule{\textwidth}{1.2pt}
~\\
\noindent [Further deformation. Hahn]. We now notice the particularly good case $4b=\lambda a^2$. It allows us to perform the next step of deformation, using the established operators as building blocks. Anyway, we set $4b=\lambda a^2$, so
$$
\varphi(y)=\frac{y}{1+\tfrac{1}{2}\lambda ay}; ~~~~~ \mathrm{U_{P^*}}=(1+\tfrac{1}{2}\lambda a\mathrm{D})^{2/\lambda} x (1+\tfrac{1}{2}\lambda a\mathrm{D})^{2-2/\lambda}=x+a(1+\lambda \theta)+\tfrac{1}{4}\lambda a^2(2+\lambda \theta)\mathrm{D}
$$
Thus our first step of deformation leads to the following dual raising and dual derivative operators
$$
\mathrm{U_{Q^*}}=\mathcal{F}_\theta^{-1}(1+\tfrac{1}{2}\lambda a\mathrm{D})^{2/\lambda} x (1+\tfrac{1}{2}\lambda a\mathrm{D})^{2-2/\lambda}\mathcal{F}_{\theta+1}; ~~~~~~~~~ \mathrm{D_{Q^*}}=\mathcal{F}_{\theta+1}^{-1} \frac{\mathrm{D}}{1-\tfrac{1}{4}\lambda^2 a^2 \mathrm{D}^2}\mathcal{F}_\theta
$$
Now we combine them as follows
\begin{align*}
\mathrm{U_{Q^*}}\mathrm{D_{Q^*}}=\mathcal{F}_\theta^{-1}(1+\tfrac{1}{2}\lambda a\mathrm{D})^{2/\lambda} &x (1+\tfrac{1}{2}\lambda a\mathrm{D})^{2-2/\lambda}\frac{\mathrm{D}}{1-\tfrac{1}{4}\lambda^2 a^2 \mathrm{D}^2}(1+\lambda\theta)\mathcal{F}_{\theta+1}=\theta_{\mathrm{Q^*}}\\
~\\
\theta_{\mathrm{Q^*}}\mathrm{U_{Q^*}}=\mathcal{F}_\theta^{-1}(1+\tfrac{1}{2}\lambda a\mathrm{D})^{2/\lambda} &x (1+\tfrac{1}{2}\lambda a\mathrm{D})^{2-2/\lambda}\frac{\mathrm{D}}{1-\tfrac{1}{4}\lambda^2 a^2 \mathrm{D}^2}(x+a(1+\lambda\theta)+\tfrac{1}{4}\lambda a^2(2+\lambda \theta)\mathrm{D})\mathcal{F}_{\theta+1}
\end{align*}
Hence we have the following almost symmetric expression
\begin{align*}
\theta_{\mathrm{Q^*}}&\mathrm{U_{Q^*}}-2a\theta_{\mathrm{Q^*}}=\\
&=\mathcal{F}_\theta^{-1}(1+\tfrac{1}{2}\lambda a\mathrm{D})^{2/\lambda} x (1+\tfrac{1}{2}\lambda a\mathrm{D})^{2-2/\lambda}\frac{\mathrm{D}}{1-\tfrac{1}{4}\lambda^2 a^2 \mathrm{D}^2}(1-\tfrac{1}{2}\lambda a\mathrm{D})^{2/\lambda} x (1-\tfrac{1}{2}\lambda a\mathrm{D})^{2-2/\lambda}\mathcal{F}_{\theta+1}
\end{align*}
and we see that the quantity in the denominator cancels. Expanding the brackets we thus obtain the following relation
$$
\theta_{\mathrm{Q^*}}\mathrm{U_{Q^*}}-2a\theta_{\mathrm{Q^*}}=\mathcal{F}_\theta^{-1}(x(1+\theta)+a(1+\lambda\theta)-a^2(1+\tfrac{1}{2}\lambda\theta)^2\mathrm{D})\mathcal{F}_{\theta+1}
$$
We may thus consider the following linear combination of established operators (still three-term)
\begin{align*}
s\mathrm{U_{Q^*}}-&\theta_{\mathrm{Q^*}}\mathrm{U_{Q^*}}+2a\theta_{\mathrm{Q^*}}=\\
&=\mathcal{F}_\theta^{-1}(x(s-1-\theta)+(s-1)a(1+\lambda\theta)+\left(\tfrac{1}{2}s\lambda a^2(1+\tfrac{1}{2}\lambda\theta)+a^2(1+\tfrac{1}{2}\lambda\theta)^2\right)\mathrm{D})\mathcal{F}_{\theta+1}
\end{align*}
Now notice that
\begin{align*}
s\mathrm{U_{Q^*}}-\theta_{\mathrm{Q^*}}\mathrm{U_{Q^*}}&+2a\theta_{\mathrm{Q^*}}=\mathrm{G_Q^{-1}} [sx-\theta x+2a\theta]\mathrm{G_Q}=\\
&=\mathrm{G_Q^{-1}}(s-1)_\theta x(1+2a\mathrm{D}) (s-1)_\theta^{-1}\mathrm{G_Q}=\mathrm{G_Q^{-1}}(s-1)_\theta \mathrm{C}_{\Delta_{2a}}^{-1} x \mathrm{C}_{\Delta_{2a}}^{\vphantom{-1}}(s-1)_\theta^{-1}\mathrm{G_Q}
\end{align*}
(Where we write $(s-1)_n=(s-1)...(s-n)$ for falling factorial, as usual and $\Delta_{2a}=(e^{2ay}-1)/2a$). Thus we may define the new family $\mathrm{S}$
$$
\mathrm{G_S}\coloneqq\mathrm{C}_{\Delta_{2a}}(s-1)_\theta^{-1}\mathrm{G_Q}(s-1)_\theta = \mathrm{C}_{\Delta_{2a}}(s-1)_\theta^{-1}\mathcal{F}_\theta^{-1}\mathrm{C}_{\mathfrak{T}f}^{-1}f'(\mathrm{D})^{-1/\lambda}\mathrm{C}_f\mathcal{F}_\theta (s-1)_\theta 
$$
that has the following three-term dual raising operator
$$
\mathrm{U_{S^*}}=x+(s-1)a+\tfrac{1}{4}a^2\frac{(2+\lambda\theta)(2+\lambda(s+\theta))(s-1-\theta)}{(1+\lambda\theta)(1+\lambda(\theta+1))}\mathrm{D}
$$
In particular case $\lambda=2, a=\tfrac{1}{2}$ we have the following deformation of Legendre polynomials, which I met in \cite{Carlitz}, defined by
$$
f^{\mathrm{S}}_0(x)=\frac{1}{s}\frac{e^{sx}-1}{e^x-1} ~~~~\iff~~~~ \mathrm{U_{S^*}}= x + \frac{s-1}{2}+\mathrm{D}\frac{\theta(s^2-\theta^2)}{4(4\theta^2-1)}
$$
These are the hypergeometric orthogonal polynomials. And the form of this expression tells us what is the expansion of these polynomials in terms of ultraspherical polynomials
\begin{align*}
p_n^{\mathrm{Q}}(x)=\mathrm{G_Q}\cdot  x^n = &\sum_{k=0}^{n} \xi_k^n x^k ~~~\Rightarrow~~~ p_n^{\mathrm{S}}(x)=\mathrm{G_S}\cdot  x^n = \mathrm{C}_{\Delta_{2a}}(s-1)_\theta^{-1}\mathrm{G_Q}(s-1)_\theta \cdot x^n =\\
&= (s-1)_n \sum_{k=0}^{n} \xi_k^n  \mathrm{C}_{\Delta_{2a}}(s-1)_\theta^{-1}\cdot x^k=\sum_{k=0}^{n} \xi_k^n  \frac{(s-1)_n}{(s-1)_k}(2a)^k \left(\frac{x}{2a}\right)_k
\end{align*}
Later we obtain further generalization of these polynomials. We established this particular case exclusively because the corresponding moment generating function is of particular interest.
~\\
~\\
\noindent\rule{\textwidth}{1.2pt}
\newpage
\noindent [Jacobi family]. We now notice that the case $4b=\lambda a^2$ actually leads to a more general family at the first step of deformation. Recall that we have fixed $f$: $f'(y)=(1+\tfrac{1}{2}\lambda af(y))^2$, $\varphi=f^{inv} \Rightarrow (\varphi')^{-1}=(1+\tfrac{1}{2}\lambda ay)^2$. Set $\mathrm{G_{P}}=f'(\mathrm{D})^{-1/\lambda}\mathrm{C}_f$, i.e. $\mathrm{G_{P}^{-1}}=(1+\tfrac{1}{2}\lambda a\mathrm{D})^{2/\lambda}\mathrm{C}_\varphi$. So that
$$
\varphi(y)=\frac{y}{1+\tfrac{1}{2}\lambda ay}; ~~~~~ \mathrm{U_{P^*}}=(1+\tfrac{1}{2}\lambda a\mathrm{D})^{2/\lambda} x (1+\tfrac{1}{2}\lambda a\mathrm{D})^{2-2/\lambda}=x+a(1+\lambda \theta)+\tfrac{1}{4}\lambda a^2(2+\lambda \theta)\mathrm{D}
$$
We are now going to use our trick again:
\begin{align*}
\mathrm{C}_{\mathfrak{T}f} x \frac{1}{1+\sigma \theta}\mathrm{C}_{\mathfrak{T}f}^{-1} = x f'(\mathrm{D})^{-1/\sigma}\frac{1}{1+\sigma x \mathfrak{T}f(\mathrm{D})}f'(\mathrm{D})^{1/\sigma}=xf'(\mathrm{D})^{-1/\sigma} \mathrm{C}_f \frac{1}{1+\sigma \theta}\mathrm{C}^{-1}_f f'(\mathrm{D})^{1/\sigma}
\end{align*}
Hence for our fixed $f$ and for any nonzero $\sigma$ the following holds
\begin{align*}
\mathrm{C}^{-1}_f f'(\mathrm{D})^{1/\lambda} \mathrm{C}_{\mathfrak{T}f} ~x \frac{1}{1+\sigma \theta}~\mathrm{C}_{\mathfrak{T}f}^{-1} f'(\mathrm{D})^{-1/\lambda} \mathrm{C}_f=\varphi'(\mathrm{D})^{-\tfrac{1}{\lambda}} ~x \varphi'(\mathrm{D})^{-1+\tfrac{1}{\sigma}}\frac{1}{1+\sigma \theta} \varphi'(\mathrm{D})^{\tfrac{1}{\lambda}-\tfrac{1}{\sigma}}
\end{align*}
We now naturally choose such values of $\sigma$ that the resulting expression contains at most three terms. The obvious one already considered is $\sigma=\lambda$, and we also have the additional one $\sigma=\kappa:~ \tfrac{1}{\kappa}-\tfrac{1}{\lambda}=\tfrac{1}{2}$. The resulting expressions are
\begin{align*}
\varphi'(\mathrm{D})^{-\tfrac{1}{\lambda}} ~x \varphi'(\mathrm{D})^{-1+\tfrac{1}{\lambda}}\frac{1}{1+\lambda \theta} &= (1+\tfrac{1}{2}\lambda a\mathrm{D})^{2/\lambda} x (1+\tfrac{1}{2}\lambda a\mathrm{D})^{2-2/\lambda} \frac{1}{1+\lambda \theta}\\
\varphi'(\mathrm{D})^{-\tfrac{1}{\lambda}} ~x \varphi'(\mathrm{D})^{-\tfrac{1}{2}+\tfrac{1}{\lambda}}\frac{1}{1+\kappa \theta} \varphi'(\mathrm{D})^{-\tfrac{1}{2}} &=  (1+\tfrac{1}{2}\lambda a\mathrm{D})^{2/\lambda} x (1+\tfrac{1}{2}\lambda a\mathrm{D})^{1-2/\lambda} \frac{1}{1+\kappa \theta} (1+\tfrac{1}{2}\lambda a\mathrm{D})
\end{align*}
or in the expanded form:
\begin{align*}
x \frac{1}{1+\lambda\theta}+a+\tfrac{1}{4}\lambda a^2 \frac{2+\lambda\theta}{1+\lambda(\theta+1)}\mathrm{D}&\\
&\text{\textit{\&}}\\
x \frac{1}{1+\kappa\theta}& + a \frac{1 + \kappa (2\theta-1)+\lambda\kappa\theta^2}{(1+\kappa(\theta-1))(1+\kappa\theta)}  +\tfrac{1}{4} \lambda a^2 \frac{2+\lambda \theta}{1+\kappa\theta}\mathrm{D}
\end{align*}
We now notice that any linear combination of two ''three-term relations'' is again a ''three-term relation''. Hence for any $r$ the following expression is still the three-term one:
\begin{align*}
\mathrm{C}^{-1}_f f'(\mathrm{D})^{1/\lambda}\mathrm{C}_{\mathfrak{T}f} ~x \left[\frac{r}{1+\lambda \theta}+\frac{1-r}{1+\kappa \theta}\right]\mathrm{C}_{\mathfrak{T}f}^{-1}f'(\mathrm{D})^{-1/\lambda} \mathrm{C}_f
\end{align*}
precisely
\begin{align*}
x\frac{1+(r\kappa +(1-r)\lambda)\theta}{(1+\lambda \theta)(1+\kappa\theta)}+~~~~~~~~&\\
+ra+(1-r)a &\frac{1 + \kappa (2\theta-1)+\lambda\kappa\theta^2}{(1+\kappa(\theta-1))(1+\kappa\theta)}+\\
&~~~~~~~~~~~~~~~~~+\tfrac{1}{4}\lambda a^2 \frac{(2+\lambda\theta)(r(1+\kappa\theta)+(1-r)(1+\lambda(\theta+1)))}{(1+\lambda (\theta+1))(1+\kappa\theta)}\mathrm{D}
\end{align*}
Now set $r\kappa+(1-r)\lambda=\beta$ and consider $\mathcal{F}_\theta = (4/(\lambda\kappa))^{\theta}\Gamma(\nicefrac{2}{\lambda})\Gamma(2\theta+\nicefrac{2}{\lambda})^{-1}\beta^{\theta}\Gamma(\theta+\nicefrac{1}{\beta})\Gamma(\nicefrac{1}{\beta})^{-1}\Rightarrow$ $\Rightarrow ~~\mathcal{F}_{\theta+1}/\mathcal{F}_\theta = (1+\beta\theta)(1+\lambda \theta)^{-1}(1+\kappa\theta)^{-1}$ (since $\tfrac{1}{\kappa}-\tfrac{1}{\lambda}=\tfrac{1}{2}$). Hence for the new family $\mathrm{Q}$
$$
\mathrm{G_Q}=\mathcal{F}^{-1}_\theta \mathrm{C}_{\mathfrak{T}f}^{-1}f'(\mathrm{D})^{-1/\lambda} \mathrm{C}_f \mathcal{F}_\theta
$$
the calculations above show that it is again orthogonal. Moreover we can recover its moment generating function as $[y^n]f_0^\mathrm{Q}(y)=\mathcal{F}_n [y^n] f'(\omega(y))^{1/\lambda}$ for $\omega=(\mathfrak{T}f)^{inv}$, hence
$$
f_0^{\mathrm{Q}}(y)=\sum_{n=0}^\infty \frac{\mathcal{F}_n}{1+\lambda n} \binom{2\lambda^{-1}+2n}{n}(\tfrac{1}{2}\lambda a y)^n=\sum_{n=0}^\infty \frac{\Gamma(\nicefrac{2}{\kappa})}{\Gamma(\nicefrac{2}{\kappa}+n)}\frac{\Gamma(\nicefrac{1}{\beta}+n)}{\Gamma(\nicefrac{1}{\beta})}\frac{(2\beta ay)^n}{\kappa^n n!}
$$
so the resulting polynomials are the classical Jacobi orthogonal polynomials $p^{\textit{J}}_n(x)=\mathrm{G_Q}\cdot x^n$.
~\\
~\\
\noindent\rule{\textwidth}{1.2pt}
\newpage
\noindent [Wilson family]. Once we established the explicit expression for $\mathrm{G_Q}$ in the case of Jacobi polynomials, it becomes a routine computation to derive explicitly the differential equation it satisfies. Consider the operator
$$
(1+\lambda\theta_\mathrm{Q})^2 = \mathrm{G_Q} (1+\lambda\theta)^2 \mathrm{G_Q^{-1}}
$$
which by definition acts by the rule $(1+\lambda\theta_\mathrm{Q})^2 \cdot p_n^{J}(x)=(1+\lambda n)^2 p_n^{J}(x)$. Now $f(y)=y(1-\tfrac{1}{2}\lambda a y)^{-1}$ hence $\mathfrak{T}f(y)=y-\tfrac{1}{2}\lambda a y^2$ and this series has the functional inverse equal to $\omega(y) = \frac{1-\sqrt{1-2\lambda a y}}{\lambda a}$, or in other words $\omega'(y)^{-2}=1-2\lambda a y$. Hence
\begin{align*}
(1+\lambda\theta_\mathrm{Q})^2 &= \mathcal{F}^{-1}_\theta \mathrm{C}_{\mathfrak{T}f}^{-1}f'(\mathrm{D})^{-1/\lambda} \mathrm{C}_f ~ (1+\lambda\theta)^2 ~ \mathrm{C}_f^{-1} f'(\mathrm{D})^{1/\lambda} \mathrm{C}_{\mathfrak{T}f}\mathcal{F}_\theta=\\
&=\mathcal{F}^{-1}_\theta (1+\lambda\theta)\omega'(\mathrm{D})^{-1}(1+\lambda\theta)\omega'(\mathrm{D})^{-1} \mathcal{F}_\theta=\\
&=\mathcal{F}^{-1}_\theta (1+\lambda\theta)((1+\lambda\theta)\omega'(\mathrm{D})^{-2}+\lambda\mathrm{D}(\omega'(\mathrm{D})^{-1})'\omega'(\mathrm{D})^{-1}) \mathcal{F}_\theta=\\
&=\mathcal{F}^{-1}_\theta (1+\lambda\theta)((1+\lambda\theta)(1-2\lambda a \mathrm{D})+\lambda\mathrm{D}(-\lambda a)) \mathcal{F}_\theta=\\
&=\mathcal{F}^{-1}_\theta [(1+\lambda\theta)^2-\tfrac{2a\lambda^2}{\kappa}(1+\lambda\theta)(1+\kappa\theta)\mathrm{D}] \mathcal{F}_\theta=\\
&=(1+\lambda\theta)^2-\tfrac{2a\lambda^2}{\kappa}(1+\beta\theta)\mathrm{D}
\end{align*}
(so that indeed $p_n^{J}(x)$ satisfies the classical hypergeometric differential equation). We now notice that since we work with conjugation, we may actually multiply any given three-term operator with this one (or take an arbitrary linear combination with this one) and the resulting operator will remain the three-term one (since it is being multiplied with $(1+\lambda\theta)^2$ or combined in any other manner). Notice that we may multiply with the first power of $(1+\lambda\theta)$ not its square, but it is difficult to deal with $\omega'(\mathrm{D})$ not being squared. We may thus deform the known Jacobi family as follows. Consider the operator
$$
(1+\lambda\theta_\mathrm{Q})^2 ~ x = (1+\lambda\theta)^2 x -\tfrac{2a\lambda^2}{\kappa}(1+\beta\theta)(1+\theta)=\mathrm{G_Q} (1+\lambda\theta)^2 \mathrm{U_Q^{*}} \mathrm{G_Q^{-1}}
$$
Obviously this is a three-term operator conjugated with $\mathrm{G_Q}$, simply the known Jacobi operator multiplied with $(1+\lambda\theta)^2$ from the left. In the expanded form we have that the following operator is three-term.
\begin{align*}
\mathrm{C}^{-1}_f f'(\mathrm{D})^{1/\lambda}\mathrm{C}_{\mathfrak{T}f} ~&\left[(1+\lambda\theta)^2 -\frac{2a\lambda^2}{\kappa}(1+\lambda\theta)(1+\kappa\theta)\mathrm{D}\right] x \left[\frac{r}{1+\lambda \theta}+\frac{1-r}{1+\kappa \theta}\right]\mathrm{C}_{\mathfrak{T}f}^{-1}f'(\mathrm{D})^{-1/\lambda} \mathrm{C}_f=\\
=\mathrm{C}^{-1}_f f'(\mathrm{D})^{1/\lambda}&\mathrm{C}_{\mathfrak{T}f} ~ \left[x\frac{(1+\beta\theta)(1+\lambda(\theta+1))^2}{(1+\lambda \theta)(1+\kappa \theta)} -\frac{2a\lambda^2}{\kappa}(1+\theta)(1+\beta\theta)\right]\mathrm{C}_{\mathfrak{T}f}^{-1}f'(\mathrm{D})^{-1/\lambda} \mathrm{C}_f
\end{align*}
Now notice that this operation allows us to consider another Jacobi family and merge it with this one. We leave $\lambda$ and $\kappa$ untouched but we choose another their combination with $\widetilde{\beta}=\widetilde{r}\kappa+(1-\widetilde{r})\lambda$. We thus have the following three-term operator
$$
\mathrm{C}^{-1}_f f'(\mathrm{D})^{1/\lambda}\mathrm{C}_{\mathfrak{T}f}\left[ [ (1+\lambda\theta)\omega'(\mathrm{D})^{-1} ]^2 x\frac{ h(1+\beta\theta)}{(1+\lambda \theta)(1+\kappa \theta)} + x\frac{(1-h) (1+\widetilde{\beta}\theta)}{(1+\lambda \theta)(1+\kappa \theta)} \right] \mathrm{C}_{\mathfrak{T}f}^{-1}f'(\mathrm{D})^{-1/\lambda} \mathrm{C}_f
$$
(affine combination of multiplied Jacobi with parameter $\beta$ and ordinary Jacobi with parameter $\widetilde{\beta}$). We now restrict our attention to the operator in brackets. Define $\mathcal{H}_\theta$: $\mathcal{H}_{\theta+1}/\mathcal{H}_\theta=(h(1+\beta\theta)(1+\lambda(\theta+1))^2+(1-h)(1+\widetilde{\beta}\theta))(1+\lambda \theta)^{-1}(1+\kappa \theta)^{-1}$ and $\tilde{\mathrm{C}}_2$: $\tilde{\mathrm{C}}_2^{-1}x \tilde{\mathrm{C}}_2 = x - 2ah\lambda^2\kappa^{-1}(1+\theta)(1+\beta\theta)$. So that the expression in brackets in the expanded form is equal to
$$
\mathcal{H}_\theta\tilde{\mathrm{C}}_2^{-1} x \tilde{\mathrm{C}}_2 \mathcal{H}_\theta^{-1} = x\frac{h(1+\beta\theta)(1+\lambda(\theta+1))^2+(1-h)(1+\widetilde{\beta}\theta)}{(1+\lambda \theta)(1+\kappa \theta)} -\frac{2ah\lambda^2}{\kappa}(1+\theta)(1+\beta\theta)
$$
Hence as we see two different Jacobi families lead to the orthogonal family
$$
\mathrm{G}_{\mathcal{W}}=\tilde{\mathrm{C}}_2 \mathcal{H}_\theta^{-1} \mathrm{C}_{\mathfrak{T}f}^{-1}f'(\mathrm{D})^{-1/\lambda} \mathrm{C}_f \mathcal{H}_\theta 
$$
The structure of this family is well understood from the form of $\mathrm{G}_{\mathcal{W}}$: we know how $\mathcal{H}_\theta$ acts and the relation $\tilde{\mathrm{C}}_\ell^{-1}x \tilde{\mathrm{C}}_\ell = x - \ell_\theta $ implies $\tilde{\mathrm{C}}_\ell \cdot x^n = \prod_{k=0}^{n-1} (x+\ell_k)$ where in the case of $\tilde{\mathrm{C}}_2$, $\ell_k$ is a quadratic polynomial in $k$. We do not symmetrize the expression above because what we really need is the feeling of evolution, how simple families create more complicated ones.
~\\
~\\
\noindent\rule{\textwidth}{1.2pt}
\newpage
\noindent [Comment on already established polynomials]. To summarize, what was actually done so far is that in addition to Jacobi generators (the first two) we have also formally obtained the following ones (the second two)
$$
x\frac{1}{1+\lambda\theta} ~~~~~~~~~~~~ x\frac{1}{1+\kappa\theta}~~~~~~~~~~~~x-2\lambda a\theta ~~~~~~~~~~~~ x(1+\lambda\theta)-a\lambda(2-\lambda+2\lambda\theta)\theta
$$
All of these operators become three-term after conjugation with $\mathrm{C}_f^{-1} f'(\mathrm{D})^{1/\lambda}\mathrm{C}_{\mathfrak{T}f}$ for our specific choice of $f$. Their affine combinations are also three-term and moreover their diagonalization is well-understood. Of course, we also have additional generators, for example
$$
x\frac{1}{\omega'(\mathrm{D})} ~~~~~~~~ x\frac{1}{\omega'(\mathrm{D})}\frac{1}{1+\lambda\theta} ~~~~~~~~ x\frac{1}{1+\kappa\theta}\frac{1}{\omega'(\mathrm{D})} ~~~~~~~~~ f'(\omega(\mathrm{D}))^{\tfrac{1}{\sigma}-\tfrac{1}{\lambda}}x\frac{1}{1+\sigma\theta} f'(\omega(\mathrm{D}))^{\tfrac{1}{\lambda}-\tfrac{1}{\sigma}}
$$
(the last one for any nonzero $\sigma$). But only in a few cases we may describe how the diagonalization of their linear combination looks like explicitly. Nevertheless there actually exist well-understood cases with additional generators, corresponding to the case $\sigma=1$ with minor changes.
\noindent\rule{\textwidth}{1.2pt}
~\\
\noindent [Long division lemma]. To proceed further, we need to establish the following fundamental property of the $0$-derivative $\mathrm{L}$. Suppose we have two power series at our disposal, say $A(x), B(x) \in 1+x\mathbb{C}[[x]]$. In addition consider two invertible operators $H_\theta$ and $R_\theta$. Suppose now that we want to understand what is the expansion of the series
$$
\frac{1}{\theta!} \cdot \frac{R_\theta \cdot A(x)}{H_\theta \cdot B(x)}
$$
We then write tautologically
\begin{align*}
\frac{1}{\theta!} \cdot \frac{R_\theta \cdot A(x)}{H_\theta \cdot B(x)}&=\exp(x \partial_y)\cdot \frac{1}{\theta_y!} \cdot \frac{R_{\theta_y} \cdot A(y)}{H_{\theta_y} \cdot B(y)}\bigg|_{y=0} = \theta_y !\exp(x \partial_y)\cdot \frac{1}{\theta_y!} \cdot \frac{R_{\theta_y} \cdot A(y)}{H_{\theta_y} \cdot B(y)}\bigg|_{y=0}=\\
&=\exp(x \mathrm{L}_y) \cdot \frac{R_{\theta_y} \cdot A(y)}{H_{\theta_y} \cdot B(y)}\bigg|_{y=0} = (H_\theta\cdot B)(y)\exp(x \mathrm{L}_y) \cdot \frac{R_{\theta_y} \cdot A(y)}{H_{\theta_y} \cdot B(y)}\bigg|_{y=0}=\\
&=\exp(x (H_\theta \cdot B)\mathrm{L}(H_\theta \cdot B)^{-1}) R_\theta \cdot A(y)\bigg|_{y=0}
\end{align*}
Now we notice that the operator $(H_\theta \cdot B)\mathrm{L}(H_\theta \cdot B)^{-1}$ in terms of $\delta=1-y\mathrm{L}$ acts as follows
$$
(H_\theta \cdot B)\mathrm{L}(H_\theta \cdot B)^{-1} \cdot f(y)=\frac{f(y)-f(0)(H_\theta \cdot B)(y)}{y} = (\mathrm{L}-\mathrm{L}(H_\theta \cdot B)\delta) \cdot f(y) = (\mathrm{L}-(\mathrm{L}H_\theta \cdot B)\delta) \cdot f(y)
$$
So we see that
$$
\frac{1}{\theta!} \cdot \frac{R_\theta \cdot A(x)}{H_\theta \cdot B(x)}=\exp\left(x\left(\dfrac{R_{\theta+1}}{R_\theta}\mathrm{L}-R_\theta^{-1}\mathrm{L}H_\theta B(y)\delta\right)\right) \cdot A(y) \bigg|_{y=0}
$$
Now lets look at this operator closely. We have $\mathrm{L}H_\theta=H_{\theta+1}\mathrm{L}$, $\delta \ell(y) = \ell(0) \delta$ and $\delta m_\theta=m_0\delta$ (so that we may tautologically multiply $\delta$ from the right by any invertible operator that has a value $1$ at $x^0$). Thus we have the following operational identity
$$
H_\theta^{-1}(H_\theta \cdot B)\mathrm{L}(H_\theta \cdot B)^{-1} H_\theta = \frac{H_{\theta+1}}{H_\theta}B\mathrm{L}B^{-1}
$$
In other words we obtained an explicit diagonalization of the operator
$$
\frac{H_{\theta+1}}{H_\theta}\mathrm{L} - t\ell(y)\delta = H_\theta^{-1} (\mathrm{L} - t(H_\theta\cdot\ell)(y)\delta) H_\theta =  H_\theta^{-1} (1+ty(H_\theta\cdot\ell))\mathrm{L}(1+ty(H_\theta\cdot\ell)(y))^{-1} H_\theta
$$
And that is very good actually, since any linear nondegenerate combination of these operators has the same form, thus also has an explicit diagonalization. The most important case for our study is the hypergeometric case, when $H_{\theta+1}/H_\theta$ is a rational function of $\theta$, that does not have any zeroes or poles at nonnegative integers. But in any case this diagonalization is preserved in polynomial domain:
$$
\ell(\mathrm{D})^{-1}x\frac{1}{1+\theta}\ell(\mathrm{D})\frac{H_{\theta+1}}{H_\theta} = H_\theta (H_\theta\cdot\ell)(\mathrm{D})^{-1}x\frac{1}{1+\theta}(H_\theta\cdot\ell)(\mathrm{D}) H_\theta^{-1}
$$
\textit{Remark}: Notice the important 'change of variable' property
$$
\mathrm{C}_f x \mathrm{C}_f^{-1}= x \frac{1}{f'(\mathrm{D})} ~~~~~~~~\text{but}~~~~~~~\mathrm{C}_f x \frac{1}{1+\theta} \mathrm{C}_f^{-1}= x \frac{1}{1+\theta}\frac{\mathrm{D}}{f(\mathrm{D})}
$$
\noindent\rule{\textwidth}{1.2pt}
~\\
\noindent [Sheffer associated family]. We are now ready to find an explicit representation of the operator $\mathrm{G_P}(c)$ of the simplest Sheffer family. Consider the associated dual raising operator.
$$
\mathrm{U_{P^*}}=x+a(1+\lambda \theta)+b(2+\lambda \theta)\mathrm{D} ~~\rightarrow~~ \mathrm{U}_{\mathrm{P}^*}(c)=x+a(1+\lambda(\theta+c))+b\frac{(2+\lambda (\theta+c))(1+\theta+c)}{1+\theta}\mathrm{D}
$$
Then conjugation with falling factorials results in the following expression
\begin{align*}
\frac{(c+\theta)_\theta}{\theta!}\mathrm{U}_{\mathrm{P}^*}(c)\frac{\theta!}{(c+\theta)_\theta}&= x\frac{(1+\theta+c)}{1+\theta}+a(1+\lambda(\theta+c))+b(2+\lambda (\theta+c))\mathrm{D}=\\
&=x+a(1+\lambda\theta)+b(2+\lambda\theta)\mathrm{D}+cx\frac{1}{1+\theta}+a\lambda c + b\lambda c\mathrm{D}=\\
&=\varphi'(\mathrm{D})^{-1/\lambda}x\varphi'(\mathrm{D})^{1/\lambda -1}+c\varphi'(\mathrm{D})^{-1}x\frac{1}{1+\theta}
\end{align*}
 (since $\mathrm{D}x=1+\theta$, and we again fix $\varphi'(y)^{-1}=1+\lambda ay+ \lambda by^2$). Further computation gives
\begin{align*}
\mathrm{C}_f\frac{(c+\theta)_\theta}{\theta!}\mathrm{U}_{\mathrm{P}^*}(c)\frac{\theta!}{(c+\theta)_\theta} &\mathrm{C}_f^{-1} = f'(\mathrm{D})^{1/\lambda}xf'(\mathrm{D})^{-1/\lambda}+ cf'(\mathrm{D})x\frac{1}{1+\theta}\frac{\mathrm{D}}{f(\mathrm{D})}=\\
 &=f'(\mathrm{D})^{1/\lambda}xf'(\mathrm{D})^{-1/\lambda}+ c\frac{f(\mathrm{D})}{\mathrm{D}}x\frac{1}{1+\theta}\frac{\mathrm{D}}{f(\mathrm{D})}+c\frac{f'(\mathrm{D})}{f(\mathrm{D})}-c\frac{1}{\mathrm{D}}=\\
&=\left(\frac{\mathrm{D}}{f(\mathrm{D})}\right)^{-c}f'(\mathrm{D})^{1/\lambda}xf'(\mathrm{D})^{-1/\lambda}\left(\frac{\mathrm{D}}{f(\mathrm{D})}\right)^{c} + c\frac{f(\mathrm{D})}{\mathrm{D}}x\frac{1}{1+\theta}\frac{\mathrm{D}}{f(\mathrm{D})}=\\
&=\left(\frac{\mathrm{D}}{f(\mathrm{D})}\right)^{-c}f'(\mathrm{D})^{1/\lambda}(c+\theta)f'(\mathrm{D})^{-1/\lambda}\left(\frac{\mathrm{D}}{f(\mathrm{D})}\right)^{c-1}x\frac{1}{1+\theta}\frac{\mathrm{D}}{f(\mathrm{D})}
\end{align*}
Hence we have obtained the following result
\begin{align*}
\left(\frac{\mathrm{D}}{f(\mathrm{D})}\right)^{c}f'(\mathrm{D})^{-1/\lambda}\mathrm{C}_f&\frac{(c+\theta)_\theta}{\theta!}\mathrm{U}_{\mathrm{P}^*}(c)\frac{\theta!}{(c+\theta)_\theta} \mathrm{C}_f^{-1}f'(\mathrm{D})^{1/\lambda}\left(\frac{\mathrm{D}}{f(\mathrm{D})}\right)^{-c}=\\
&=(c+\theta)f'(\mathrm{D})^{-1/\lambda}\left(\frac{\mathrm{D}}{f(\mathrm{D})}\right)^{c-1}x\frac{1}{1+\theta}\left(\frac{\mathrm{D}}{f(\mathrm{D})}\right)^{1-c}f'(\mathrm{D})^{1/\lambda}
\end{align*}
We now denote $\ell(y)=f'(y)^{1/\lambda}(y/f(y))^{1-c}$ and consider the invertible operator $H_\theta = (-1)^\theta(-c)_\theta =$ $=(c+\theta-1)_\theta$, so that $H_{\theta+1}/H_\theta = c+\theta$. Then by the long division lemma we have the equality
\begin{align*}
\left(\frac{\mathrm{D}}{f(\mathrm{D})}\right)^{c}&f'(\mathrm{D})^{-1/\lambda}\mathrm{C}_f\frac{(c+\theta)_\theta}{\theta!}\mathrm{U}_{\mathrm{P}^*}(c)\frac{\theta!}{(c+\theta)_\theta} \mathrm{C}_f^{-1}f'(\mathrm{D})^{1/\lambda}\left(\frac{\mathrm{D}}{f(\mathrm{D})}\right)^{-c}=\\
&=\frac{H_{\theta+1}}{H_\theta}\ell(\mathrm{D})^{-1}x\frac{1}{1+\theta}\ell(\mathrm{D})=H_{\theta+1} (H_\theta \cdot \ell)(\mathrm{D})^{-1}x\frac{1}{1+\theta}(H_\theta \cdot \ell)(\mathrm{D})H_{\theta+1}^{-1}
\end{align*}
So that we have established the following explicit representation
\begin{align*}
\mathrm{G_P}(c)=\theta! [H_\theta \cdot\ell](\mathrm{D})(c+\theta)_\theta^{-1}\left(\frac{\mathrm{D}}{f(\mathrm{D})}\right)^{c}&f'(\mathrm{D})^{-1/\lambda}\mathrm{C}_f\frac{(c+\theta)_\theta}{\theta!}
\end{align*}
Now with this explicit representation at our disposal it becomes an easy task to compute the dual polynomials $\mathrm{G_P}(c)^{-1}\cdot x^n$ or the moment generating function
$$
\overline{\mathrm{G_P}(c)^{-1}}\cdot 1=f_0(y,c)=\frac{1}{\theta!} \cdot \frac{(c+\theta)_\theta \cdot f'(y)^{1/\lambda}(y/f(y))^{-c}}{(c+\theta-1)_\theta \cdot f'(y)^{1/\lambda}(y/f(y))^{1-c}}
$$
~\\
\noindent\rule{\textwidth}{1.2pt}
\newpage
\noindent [Ultraspherical associated family]. We now proceed in the same manner as in case of usual deformation Sheffer $\rightarrow$ Ultraspherical. I.e. we want to find an explicit form of an operator $\mathrm{G}_\mu$, such that\\
$$
\mathrm{G}_\mu^{-1} x\frac{1}{1+\lambda(\theta+c)} \mathrm{G}_\mu = x\frac{1}{1+\lambda(\theta_\mathrm{P}(c)+c)}
$$
We again write $f'(y)=1+\lambda a f(y)+\lambda b f^2(y)$ and again this property of $f$ is required only for ''three-term-ness'', not for the structure of conjugations. With an easy guess we proceed as follows
\begin{align*}
\mathrm{U}_{\mathrm{P}^*}&(c) \frac{1}{1+\lambda(\theta+c)} =\\
&=\frac{\theta!}{(c+\theta)_\theta} \left[\left[\varphi'(\mathrm{D})^{-1/\lambda}x\varphi'(\mathrm{D})^{1/\lambda -1}+c\varphi'(\mathrm{D})^{-1}x\frac{1}{1+\theta}\right]\frac{1}{1+\lambda(\theta+c)}\right] \frac{(c+\theta)_\theta}{\theta!}=\\
&=\frac{\theta!}{(c+\theta)_\theta} \mathrm{C}^{-1}_f\left[\left[f'(\mathrm{D})^{1/\lambda}xf'(\mathrm{D})^{-1/\lambda}+cf'(\mathrm{D})x\frac{1}{1+\theta}\frac{\mathrm{D}}{f(\mathrm{D})}\right]\frac{1}{1+\lambda(c+x\mathfrak{T}f(\mathrm{D}))}\right]\mathrm{C}_f  \frac{(c+\theta)_\theta}{\theta!}
\end{align*}
Hence the additional conjugation results in the following equality.
\begin{align*}
f'(\mathrm{D})&^{-c-\tfrac{1}{\lambda}} \mathrm{C}_f \frac{(c+\theta)_\theta}{\theta!}\left[\mathrm{U}_{\mathrm{P}^*}(c) \frac{1}{1+\lambda(\theta+c)}\right]\frac{\theta!}{(c+\theta)_\theta} \mathrm{C}^{-1}_f f'(\mathrm{D})^{c+\tfrac{1}{\lambda}}=\\
&=\left[f'(\mathrm{D})^{-c}xf'(\mathrm{D})^{c}+cf'(\mathrm{D})^{1-\tfrac{1}{\lambda}-c}x\frac{1}{1+\theta}\frac{\mathrm{D}}{f(\mathrm{D})} f'(\mathrm{D})^{c+\tfrac{1}{\lambda}}\right]\frac{1}{(1+\lambda c)(\mathfrak{T}f)'(\mathrm{D})+x\mathfrak{T}f(\mathrm{D})}
\end{align*}
We now again introduce $\omega=(f/f')^{inv}$. Now we have $(1-(f/f')'(y))/(f/f')=f''(y)/f'(y)$ and hence $(f'(\omega)^c)'f'(\omega)^{-c}= c(\omega'-1)/y$. Thus the following holds
\begin{align*}
\mathrm{C}_{\mathfrak{T}f}^{-1}&f'(\mathrm{D})^{-c-\tfrac{1}{\lambda}}\mathrm{C}_f \frac{(c+\theta)_\theta}{\theta!}\left[\mathrm{U}_{\mathrm{P}^*}(c) \frac{1}{1+\lambda(\theta+c)}\right]\frac{\theta!}{(c+\theta)_\theta} \mathrm{C}^{-1}_f f'(\mathrm{D})^{c+\tfrac{1}{\lambda}}\mathrm{C}_{\mathfrak{T}f}=\\
&=\left[f'(\omega(\mathrm{D}))^{-c}xf'(\omega(\mathrm{D}))^{c}+cf'(\omega(\mathrm{D}))^{1-\tfrac{1}{\lambda}-c}x\frac{1}{1+\theta} f'(\omega(\mathrm{D}))^{c+\tfrac{1}{\lambda}-1}\omega'(\mathrm{D})\right]\frac{1}{1+\lambda(\theta+c)}=\\
&=\left[x-c\frac{\omega'(\mathrm{D})-1}{\mathrm{D}}+cf'(\omega(\mathrm{D}))^{1-\tfrac{1}{\lambda}-c}x\frac{1}{1+\theta} f'(\omega(\mathrm{D}))^{c+\tfrac{1}{\lambda}-1}\omega'(\mathrm{D})\right]\frac{1}{1+\lambda(\theta+c)}=\\
&=\left[x+c\omega'(\mathrm{D})^{-1}f'(\omega(\mathrm{D}))^{1-\tfrac{1}{\lambda}-c}x\frac{1}{1+\theta} f'(\omega(\mathrm{D}))^{c+\tfrac{1}{\lambda}-1}\omega'(\mathrm{D})\right]\frac{1}{1+\lambda(\theta+c)}=\\
&=  (c+\theta)\left[\omega'(\mathrm{D})f'(\omega(\mathrm{D}))^{c+\tfrac{1}{\lambda}-1}\right]^{-1} x\frac{1}{1+\theta} \left[\omega'(\mathrm{D})f'(\omega(\mathrm{D}))^{c+\tfrac{1}{\lambda}-1}\right] \frac{1}{1+\lambda(\theta+c)}
\end{align*}
We now denote $\ell(y)=\omega'(y)f'(\omega(y))^{c+\tfrac{1}{\lambda}-1}$, and introduce two sequences $\mathcal{F}_\theta = \lambda^{-\theta}\Gamma(c+\nicefrac{1}{\lambda})\Gamma(c+\theta+\nicefrac{1}{\lambda})^{-1}~\Rightarrow~ \mathcal{F}_{\theta+1}/\mathcal{F}_\theta=F_\theta=(1+\lambda (\theta+c))^{-1}$ and  $H_\theta = (-1)^\theta(-c)_\theta=(c+\theta-1)_\theta$, so that $H_{\theta+1}/H_\theta = c+\theta$. Then by the long divison lemma we have the equality
\begin{align*}
\mathrm{C}_{\mathfrak{T}f}^{-1}&f'(\mathrm{D})^{-c-\tfrac{1}{\lambda}}\mathrm{C}_f \frac{(c+\theta)_\theta}{\theta!}\left[\mathrm{U}_{\mathrm{P}^*}(c) \frac{1}{1+\lambda(\theta+c)}\right]\frac{\theta!}{(c+\theta)_\theta} \mathrm{C}^{-1}_f f'(\mathrm{D})^{c+\tfrac{1}{\lambda}}\mathrm{C}_{\mathfrak{T}f}=\\
&=\frac{H_{\theta+1}}{H_\theta}\ell(\mathrm{D})^{-1}x\frac{1}{1+\theta}\ell(\mathrm{D})\frac{\mathcal{F}_{\theta+1}}{\mathcal{F}_\theta}=H_{\theta+1} (H_\theta \cdot \ell)(\mathrm{D})^{-1}x\frac{1}{1+\theta}(H_\theta \cdot \ell)(\mathrm{D})H_{\theta+1}^{-1}\frac{\mathcal{F}_{\theta+1}}{\mathcal{F}_\theta}=\\
&=\mathcal{F}_\theta H_{\theta+1} (H_\theta\mathcal{F}_\theta \cdot \ell)(\mathrm{D})^{-1}x\frac{1}{1+\theta}(H_\theta\mathcal{F}_\theta\cdot \ell)(\mathrm{D})H_{\theta+1}^{-1} \mathcal{F}_\theta^{-1}
\end{align*}
So that we have established the following explicit representation for ultraspherical associated polynomials
$$
\mathrm{G_Q}(c)=\theta! (H_\theta\mathcal{F}_\theta\cdot \ell)(\mathrm{D})(c+\theta)_{\theta}^{-1} \mathcal{F}_\theta^{-1}\mathrm{C}_{\mathfrak{T}f}^{-1}f'(\mathrm{D})^{-c-\tfrac{1}{\lambda}}\mathrm{C}_f \mathcal{F}_\theta \frac{(c+\theta)_\theta}{\theta!}
$$
with the corresponding dual raising operator explicitly given by
$$
\mathrm{U_{Q^*}}(c)=\mathcal{F}_\theta^{-1} \mathrm{U_{P^*}}(c) \mathcal{F}_{\theta+1}=x+a+b\frac{(2+\lambda (\theta+c))(1+\theta+c)}{(1+\theta)(1+\lambda(\theta+c))(1+\lambda(\theta+c+1))}\mathrm{D}
$$
One may now compute the series $H_\theta\mathcal{F}_\theta\cdot \ell(y)$ by the use of combinatorial residue theorem. The corresponding moment generating function is a well-understood fraction. Clearly the case $c=0$ produces usual ultraspherical polynomials.
~\\
\noindent\rule{\textwidth}{1.2pt}
\newpage
\noindent [Simple observation]. We now make the following observation. Suppose that we have a family $\mathrm{P}$ with the dual raising operator of the form
$$
\mathrm{U_{P^*}}=x+\sum_{i=1}^{n}\alpha_i k_i (\theta)+\left[\sum_{i=1}^{n}\alpha_i s_i (\theta)\right]\left[\sum_{i=1}^{n}\alpha_i h_i (\theta)\right]\mathrm{D}
$$
and $\mathcal{K}_{\theta+1}/\mathcal{K}_\theta=\sum_{i=1}^{n}\alpha_i h_i (\theta)$ is invertible. In other words, the following operator splits
$$
\mathcal{K}_\theta \mathrm{U_{P^*}}\mathcal{K}_\theta^{-1} = \sum_{i=1}^{n}\alpha_i [xh_i(\theta)+k_i(\theta)+s_i(\theta)\mathrm{D}]
$$
Then the associated dual raising operator
$$
\mathrm{U_{P^*}}(c)=x+\sum_{i=1}^{n}\alpha_i k_i (\theta+c)+\frac{1+\theta+c}{1+\theta}\left[\sum_{i=1}^{n}\alpha_i s_i (\theta+c)\right]\left[\sum_{i=1}^{n}\alpha_i h_i (\theta+c)\right]\mathrm{D}
$$
also splits after conjugation, if $\mathcal{K}_{\theta+c}\mathcal{K}_{c}^{-1}$ is also invertible.
$$
\frac{\mathcal{K}_{\theta+c}}{\mathcal{K}_{c}}\frac{(c+\theta)_\theta}{\theta!}\mathrm{U_{P^*}}(c)\frac{\theta!}{(c+\theta)_\theta}\frac{\mathcal{K}_{c}}{\mathcal{K}_{\theta+c}} = \sum_{i=1}^{n}\alpha_i \left[x\frac{1+\theta+c}{1+\theta}h_i(\theta+c)+k_i(\theta+c)+s_i(\theta+c)\mathrm{D}\right]
$$
\noindent\rule{\textwidth}{1.2pt}
~\\
\noindent [Jacobi associated family]. We now consider the case $4b=\lambda a^2$, so that $f'(y)=(1+\tfrac{1}{2}\lambda af(y))^{2}$, $\varphi =f^{inv}=y(1+\tfrac{1}{2}\lambda a y)^{-1}$, $f/f'=y-\tfrac{1}{2}\lambda ay^2$, $\omega = (f/f')^{inv} = \tfrac{1-\sqrt{1-2\lambda a y}}{\lambda a}$. Recall that the Jacobi family is defined by
$$
\mathrm{G}_\mathcal{J}=\mathcal{F}^{-1}_\theta \mathrm{C}_{\mathfrak{T}f}^{-1}f'(\mathrm{D})^{-1/\lambda} \mathrm{C}_f \mathcal{F}_\theta
$$
where $\tfrac{1}{\kappa}-\tfrac{1}{\lambda}=\tfrac{1}{2}$ and $\mathcal{F}_{\theta+1}\mathcal{F}_\theta^{-1}=(1+\beta\theta)(1+\lambda\theta)^{-1}(1+\kappa\theta)^{-1}$, ($\mathcal{F}_0=1$). Its dual raising operator splits after conjugation with $\mathcal{F}_\theta$ into $\lambda$-part and $\kappa$-part.
\begin{align*}
\mathcal{F}_{\theta} \mathrm{U}_{\mathcal{J}^*}\mathcal{F}^{-1}_\theta = r&\left[x \frac{1}{1+\lambda\theta}+a+\tfrac{1}{4}\lambda a^2 \frac{2+\lambda\theta}{1+\lambda(\theta+1)}\mathrm{D}\right] +\\
&+(1-r)\left[x \frac{1}{1+\kappa\theta} + \frac{1}{2}a \frac{2+\lambda\theta}{1+\kappa\theta} +\frac{1}{2}a \frac{\lambda\theta}{1+\kappa(\theta-1)}+\tfrac{1}{4} \lambda a^2 \frac{2+\lambda \theta}{1+\kappa\theta}\mathrm{D}\right]
\end{align*}
Thus according to our observation, we also have the following splitting
\begin{align*}
&\frac{\mathcal{F}_{\theta+c}}{\mathcal{F}_{c}}\frac{(c+\theta)_\theta}{\theta!}\mathrm{U}_{\mathcal{J}^*}(c)\frac{\theta!}{(c+\theta)_\theta}\frac{\mathcal{F}_{c}}{\mathcal{F}_{\theta+c}}= r\left[x \frac{1+\theta+c}{(1+\theta)(1+\lambda(\theta+c))}+a+\tfrac{1}{4}\lambda a^2 \frac{2+\lambda(\theta+c)}{1+\lambda(1+\theta+c)}\mathrm{D}\right] +\\
&+(1-r)\left[x \frac{1+\theta+c}{(1+\theta)(1+\kappa(\theta+c))} + \tfrac{1}{2}a \frac{2+\lambda(\theta+c)}{1+\kappa(\theta+c)} +\tfrac{1}{2}a \frac{\lambda(\theta+c)}{1+\kappa(\theta+c-1)}+\tfrac{1}{4} \lambda a^2 \frac{2+\lambda (\theta+c)}{1+\kappa(\theta+c)}\mathrm{D}\right]
\end{align*}
Now the first summand ('$\lambda$-part') is something familiar, it corresponds to the ultraspherical case. Thus we know the operator that conjugation with $\mathrm{C}_{\mathfrak{T}f}^{-1} f'(\mathrm{D})^{-c-1/\lambda}\mathrm{C}_f $  leads to in this case. The question is, what happens in $r=0$ case, i.e. with the second operator. Consider the following computation.
\begin{align*}
&\mathrm{C}^{-1}_f f'(\mathrm{D})^{c+\tfrac{1}{\lambda}}\mathrm{C}_{\mathfrak{T}f} \left[x\frac{1}{1+\kappa(\theta+c)}\right] \mathrm{C}_{\mathfrak{T}f}^{-1} f'(\mathrm{D})^{-c-\tfrac{1}{\lambda}}\mathrm{C}_f =\\
&~~~~~~~~= \mathrm{C}^{-1}_f f'(\mathrm{D})^{c+\tfrac{1}{\lambda}}\left[x\frac{1}{(1+\kappa c)(\mathfrak{T}f)'(\mathrm{D})+\kappa x \mathfrak{T}f(\mathrm{D})}\right] f'(\mathrm{D})^{-c-\tfrac{1}{\lambda}}\mathrm{C}_f=\\
&~~~~~~~~= \mathrm{C}^{-1}_f\left[f'(\mathrm{D})^{c+\tfrac{1}{\lambda}} x f'(\mathrm{D})^{-c-\tfrac{1}{\kappa}}\frac{1}{1+\kappa (c+ x\mathfrak{T}f(\mathrm{D}))} f'(\mathrm{D})^{\tfrac{1}{2}}\right]\mathrm{C}_f=\\
&~~~~~~~~=\varphi'(\mathrm{D})^{-c-\tfrac{1}{\lambda}} x \varphi'(\mathrm{D})^{c+\tfrac{1}{\lambda}-\tfrac{1}{2}}\frac{1}{1+\kappa (\theta+c)} \varphi'(\mathrm{D})^{-\tfrac{1}{2}}=\\
&=x \frac{1}{1+\kappa(\theta+c)} + \tfrac{1}{2}a \frac{2+\lambda(\theta+2c)}{1+\kappa(\theta+c)}+\tfrac{1}{2} a \frac{\lambda\theta}{1+\kappa(c+\theta-1)}+\tfrac{1}{4}\lambda a^2 \frac{2+\lambda(\theta+2c)}{1+\kappa(\theta+c)}\mathrm{D}
\end{align*}
The result is quite close to what we need, so we subtract two operators
\begin{align*}
\frac{\mathcal{F}_{\theta+c}}{\mathcal{F}_{c}}&\frac{(c+\theta)_\theta}{\theta!}\mathrm{U}_{\mathcal{J}^*}(c)\frac{\theta!}{(c+\theta)_\theta}\frac{\mathcal{F}_{c}}{\mathcal{F}_{\theta+c}}\bigg|_{r=0}-\mathrm{C}^{-1}_f f'(\mathrm{D})^{c+\tfrac{1}{\lambda}}\mathrm{C}_{\mathfrak{T}f} \left[x\frac{1}{1+\kappa(\theta+c)}\right] \mathrm{C}_{\mathfrak{T}f}^{-1} f'(\mathrm{D})^{-c-\tfrac{1}{\lambda}}\mathrm{C}_f=\\
&= x \frac{c}{(1+\theta)(1+\kappa(\theta+c))} - \tfrac{1}{2}a \frac{\lambda c}{1+\kappa(\theta+c)}+\tfrac{1}{2} a \frac{\lambda c}{1+\kappa(c+\theta-1)}-\tfrac{1}{4}\lambda a^2 \frac{\lambda c}{1+\kappa(\theta+c)}\mathrm{D}=\\
&=c(1-\tfrac{1}{2}\lambda a\mathrm{D})\frac{1}{1+\kappa(c+\theta-1)}(1+\tfrac{1}{2}\lambda a\mathrm{D}) x \frac{1}{1+\theta}=\\
&=c(\mathfrak{T}f)'(\varphi(\mathrm{D}))\varphi'(\mathrm{D})^{-\tfrac{1}{2}}\frac{1}{1+\kappa(c+\theta-1)}\varphi'(\mathrm{D})^{-\tfrac{1}{2}}x \frac{1}{1+\theta}
\end{align*}
And it is left to apply conjugation with $\mathrm{C}_{\mathfrak{T}f}^{-1} f'(\mathrm{D})^{-c-1/\lambda}\mathrm{C}_f$ again.
\begin{align*}
&\mathrm{C}_{\mathfrak{T}f}^{-1}f'(\mathrm{D})^{-c-\tfrac{1}{\lambda}}\mathrm{C}_f \left[(\mathfrak{T}f)'(\varphi(\mathrm{D}))\varphi'(\mathrm{D})^{-\tfrac{1}{2}}\frac{1}{1+\kappa(c+\theta-1)}\varphi'(\mathrm{D})^{-\tfrac{1}{2}}x \frac{1}{1+\theta}\right] \mathrm{C}^{-1}_f f'(\mathrm{D})^{c+\tfrac{1}{\lambda}}\mathrm{C}_{\mathfrak{T}f}=\\
&=\mathrm{C}_{\mathfrak{T}f}^{-1}f'(\mathrm{D})^{-c-\tfrac{1}{\lambda}}\left[(\mathfrak{T}f)'(\mathrm{D})f'(\mathrm{D})^{\tfrac{1}{2}}\frac{1}{1+\kappa(c-1)+\kappa x\mathfrak{T}f(\mathrm{D})}f'(\mathrm{D})^{\tfrac{1}{2}}x \frac{1}{1+\theta}\frac{\mathrm{D}}{f(\mathrm{D})}\right] f'(\mathrm{D})^{c+\tfrac{1}{\lambda}}\mathrm{C}_{\mathfrak{T}f}=\\
&=\mathrm{C}_{\mathfrak{T}f}^{-1}\left[(\mathfrak{T}f)'(\mathrm{D})f'(\mathrm{D})^{1-c-\tfrac{1}{\kappa}}\frac{1}{1+\kappa(c-1)+\kappa x\mathfrak{T}f(\mathrm{D})}f'(\mathrm{D})^{\tfrac{1}{2}}x \frac{1}{1+\theta}\frac{\mathrm{D}}{f(\mathrm{D})} f'(\mathrm{D})^{c+\tfrac{1}{\lambda}}\right]\mathrm{C}_{\mathfrak{T}f}=\\
&=\mathrm{C}_{\mathfrak{T}f}^{-1}\left[(\mathfrak{T}f)'(\mathrm{D})\frac{1}{(1+\kappa(c-1))(\mathfrak{T}f)'(\mathrm{D})+\kappa x\mathfrak{T}f(\mathrm{D})}f'(\mathrm{D})^{1-c-\tfrac{1}{\lambda}}x \frac{1}{1+\theta}\frac{\mathrm{D}}{f(\mathrm{D})} f'(\mathrm{D})^{c+\tfrac{1}{\lambda}}\right]\mathrm{C}_{\mathfrak{T}f}=\\
&=\frac{1}{1+\kappa(c+\theta-1)}\left[f'(\omega(\mathrm{D}))^{c+\tfrac{1}{\lambda}-1}\right]^{-1} x\frac{1}{1+\theta} \left[f'(\omega(\mathrm{D}))^{c+\tfrac{1}{\lambda}-1}\right]
\end{align*}
So that we have established the following equality
\begin{align*}
\mathrm{C}_{\mathfrak{T}f}^{-1}f'(\mathrm{D})^{-c-\tfrac{1}{\lambda}}\mathrm{C}_f& \left[ \frac{\mathcal{F}_{\theta+c}}{\mathcal{F}_{c}}\frac{(c+\theta)_\theta}{\theta!}\mathrm{U}_{\mathcal{J}^*}(c)\frac{\theta!}{(c+\theta)_\theta}\frac{\mathcal{F}_{c}}{\mathcal{F}_{\theta+c}}\bigg|_{r=0} \right] \mathrm{C}^{-1}_f f'(\mathrm{D})^{c+\tfrac{1}{\lambda}}\mathrm{C}_{\mathfrak{T}f} = \\
&=x\frac{1}{1+\kappa(\theta+c)} + c\frac{1}{1+\kappa(c+\theta-1)}f'(\omega(\mathrm{D}))^{1-c-\tfrac{1}{\lambda}}x \frac{1}{1+\theta} f'(\omega(\mathrm{D}))^{c+\tfrac{1}{\lambda}-1}=\\
&=\frac{c+\theta}{1+\kappa(c+\theta-1)}\left[f'(\omega(\mathrm{D}))^{c+\tfrac{1}{\lambda}-1}\right]^{-1} x\frac{1}{1+\theta} \left[f'(\omega(\mathrm{D}))^{c+\tfrac{1}{\lambda}-1}\right]
\end{align*}
And we are almost done. Consider again the ultraspherical case
\begin{align*}
\mathrm{C}_{\mathfrak{T}f}^{-1}f'(\mathrm{D})^{-c-\tfrac{1}{\lambda}}\mathrm{C}_f& \left[ \frac{\mathcal{F}_{\theta+c}}{\mathcal{F}_{c}}\frac{(c+\theta)_\theta}{\theta!}\mathrm{U}_{\mathcal{J}^*}(c)\frac{\theta!}{(c+\theta)_\theta}\frac{\mathcal{F}_{c}}{\mathcal{F}_{\theta+c}}\bigg|_{r=1} \right] \mathrm{C}^{-1}_f f'(\mathrm{D})^{c+\tfrac{1}{\lambda}}\mathrm{C}_{\mathfrak{T}f} =\\
&=(c+\theta)\left[\omega'(\mathrm{D})f'(\omega(\mathrm{D}))^{c+\tfrac{1}{\lambda}-1}\right]^{-1} x\frac{1}{1+\theta} \left[\omega'(\mathrm{D})f'(\omega(\mathrm{D}))^{c+\tfrac{1}{\lambda}-1}\right] \frac{1}{1+\lambda(\theta+c)}
\end{align*}
Notice that we may simplify this expression further. We notice that $(1+\lambda (\theta+c-1))\cdot f'(\omega)^{-1+c+1/\lambda}=$ $=(1+\lambda (c-1))\omega'(\mathrm{D})f'(\omega(\mathrm{D}))^{c-1+1/\lambda}$ (again since $\omega'(f''/f')(\omega)=(\omega'-1)/y$). Now by long division lemma for any invertible $m_\theta$ we have the following transformation
$$
B\mathrm{L}B^{-1} m_\theta = m_{\theta+1} \left(\frac{m_0}{m_\theta}\cdot B\right)\mathrm{L}\left(\frac{m_0}{m_\theta}\cdot B\right)^{-1}
$$
Thus setting $B=f'(\omega)^{c-1+1/\lambda}$ and $m_\theta = (1+\lambda(\theta+c-1))^{-1}$ we obtain the following expression in ultraspherical case.
\begin{align*}
\mathrm{C}_{\mathfrak{T}f}^{-1}f'(\mathrm{D})^{-c-\tfrac{1}{\lambda}}\mathrm{C}_f& \left[ \frac{\mathcal{F}_{\theta+c}}{\mathcal{F}_{c}}\frac{(c+\theta)_\theta}{\theta!}\mathrm{U}_{\mathcal{J}^*}(c)\frac{\theta!}{(c+\theta)_\theta}\frac{\mathcal{F}_{c}}{\mathcal{F}_{\theta+c}}\bigg|_{r=1} \right] \mathrm{C}^{-1}_f f'(\mathrm{D})^{c+\tfrac{1}{\lambda}}\mathrm{C}_{\mathfrak{T}f} =\\
&= \frac{c+\theta}{1+\lambda(c+\theta-1)}\left[f'(\omega(\mathrm{D}))^{c+\tfrac{1}{\lambda}-1}\right]^{-1} x\frac{1}{1+\theta} \left[f'(\omega(\mathrm{D}))^{c+\tfrac{1}{\lambda}-1}\right]
\end{align*}
Now $\beta=r\kappa+(1-r)\lambda$. Denote also $k(y)=f'(\omega(y))^{c-1+1/\lambda}$ and $H_\theta = (-1)^\theta(-c)_\theta=(c+\theta-1)_\theta$, so that $H_{\theta+1}/H_\theta = c+\theta$. We may finally write the full expansion
\begin{align*}
\mathrm{C}_{\mathfrak{T}f}^{-1}f'(\mathrm{D}&)^{-c-\tfrac{1}{\lambda}}\mathrm{C}_f \left[ \frac{\mathcal{F}_{\theta+c}}{\mathcal{F}_{c}}\frac{(c+\theta)_\theta}{\theta!}\mathrm{U}_{\mathcal{J}^*}(c)\frac{\theta!}{(c+\theta)_\theta}\frac{\mathcal{F}_{c}}{\mathcal{F}_{\theta+c}} \right] \mathrm{C}^{-1}_f f'(\mathrm{D})^{c+\tfrac{1}{\lambda}}\mathrm{C}_{\mathfrak{T}f} =\\
&= \frac{(c+\theta)(1+\beta(c+\theta-1))}{(1+\lambda(c+\theta-1))(1+\kappa(c+\theta-1))}\left[f'(\omega(\mathrm{D}))^{c+\tfrac{1}{\lambda}-1}\right]^{-1} x\frac{1}{1+\theta} \left[f'(\omega(\mathrm{D}))^{c+\tfrac{1}{\lambda}-1}\right]=\\
&= \frac{H_{\theta+1}}{H_\theta}\frac{\mathcal{F}_{\theta+c}}{\mathcal{F}_{\theta+c-1}}k(\mathrm{D})^{-1} x\frac{1}{1+\theta}k(\mathrm{D})=\\
&=\frac{H_{\theta+1}}{H_1}\frac{\mathcal{F}_{\theta+c}}{\mathcal{F}_c} \left(H_\theta \frac{\mathcal{F}_{\theta+c-1}}{\mathcal{F}_{c-1}} \cdot k\right)^{-1}(\mathrm{D})x\frac{1}{1+\theta} \left(H_\theta \frac{\mathcal{F}_{\theta+c-1}}{\mathcal{F}_{c-1}} \cdot k\right)(\mathrm{D}) \frac{\mathcal{F}_c}{\mathcal{F}_{\theta+c}}\frac{H_1}{H_{\theta+1}}
\end{align*}
So that we have established the following explicit $\mathrm{G_P}$-operator for associated Jacobi family
$$
\mathrm{G}_\mathcal{J}(c)=\theta! \left(H_\theta \frac{\mathcal{F}_{\theta+c-1}}{\mathcal{F}_{c-1}} \cdot f'(\omega)^{c-1+1/\lambda}\right)(\mathrm{D})\frac{\mathcal{F}_c}{\mathcal{F}_{\theta+c}}(c+\theta)_\theta^{-1}\mathrm{C}_{\mathfrak{T}f}^{-1}f'(\mathrm{D})^{-c-\tfrac{1}{\lambda}}\mathrm{C}_f\frac{\mathcal{F}_{\theta+c}}{\mathcal{F}_{c}}\frac{(c+\theta)_\theta}{\theta!}
$$
$\mathrm{G}_\mathcal{J}(c) \cdot x^n$ are the corresponding orthogonal polynomials. The moment generating function is the inverse Laplace transform of the fraction of two $_2F_1$-hypergeometric series.
\begin{align*}
\theta!\cdot\overline{\mathrm{G}_\mathcal{J}(c)^{-1}}\cdot 1 =F_0^{\mathcal{J}}(y,c)=\frac{(c+\theta)_{\theta} \tfrac{\mathcal{F}_{\theta+c}}{\mathcal{F}_{c}} \cdot f'(\omega(y))^{c+1/\lambda}}{(c-1+\theta)_\theta\tfrac{\mathcal{F}_{\theta+c-1}}{\mathcal{F}_{c-1}} \cdot f'(\omega(y))^{c-1+1/\lambda}}=\frac{_2 F_1\left(\begin{matrix} c+1, c+\nicefrac{1}{\beta}\\ 2c+\nicefrac{2}{\kappa}\end{matrix}\bigg|\dfrac{2\beta ay}{\kappa}\right)}{_2 F_1\left(\begin{matrix} c , c-1+\nicefrac{1}{\beta}\\ 2c-2+\nicefrac{2}{\kappa}\end{matrix}\bigg|\dfrac{2\beta ay}{\kappa}\right)}
\end{align*}
\noindent\rule{\textwidth}{1.2pt}
~\\
\noindent [Wilson associated family]. Recall that we define Wilson polynomials by the use of linear combination of the dual raising operators of ordinary Jacobi polynomials and of Jacobi polynomials 'multiplied with differential equation they satisfy'. So that for $\beta=r\kappa+(1-r)\lambda$, $\widetilde{\beta}=\widetilde{r}\kappa+(1-\widetilde{r})\lambda$, invertible $\mathcal{H}_\theta$: $\mathcal{H}_{\theta+1}/\mathcal{H}_\theta=(h(1+\beta\theta)(1+\lambda(\theta+1))^2+(1-h)(1+\widetilde{\beta}\theta))(1+\lambda \theta)^{-1}(1+\kappa \theta)^{-1}$, ($\mathcal{H}_0=1$) and $\tilde{\mathrm{C}}_2$: $\tilde{\mathrm{C}}_2^{-1}x \tilde{\mathrm{C}}_2 = x - 2ah\lambda^2\kappa^{-1}(1+\theta)(1+\beta\theta)$ we define
$$
\mathrm{G}_{\mathcal{W}}=\tilde{\mathrm{C}}_2 \mathcal{H}_\theta^{-1} \mathrm{C}_{\mathfrak{T}f}^{-1}f'(\mathrm{D})^{-1/\lambda} \mathrm{C}_f \mathcal{H}_\theta 
$$
and we have the splitting
\begin{align*}
&\mathcal{H}_\theta\mathrm{U}_{\mathcal{W}^*}\mathcal{H}_\theta^{-1} =\\
&=\mathrm{C}^{-1}_f f'(\mathrm{D})^{1/\lambda}\mathrm{C}_{\mathfrak{T}f}\left[ [ (1+\lambda\theta)\omega'(\mathrm{D})^{-1} ]^2 x\frac{ h(1+\beta\theta)}{(1+\lambda \theta)(1+\kappa \theta)} + x\frac{(1-h) (1+\widetilde{\beta}\theta)}{(1+\lambda \theta)(1+\kappa \theta)} \right] \mathrm{C}_{\mathfrak{T}f}^{-1}f'(\mathrm{D})^{-1/\lambda} \mathrm{C}_f
\end{align*}
We now want to find the $\mathrm{G_P}$-operator of associated Wilson family. The case $h=0$ corresponds to associated Jacobi polynomials, thus it is well-understood. Hence according to our observation, it is enough to understand particular case $h=1$ to understand the whole general case. But we have
\begin{align*}
&\mathcal{H}_\theta\mathrm{U}_{\mathcal{W}^*}\mathcal{H}_\theta^{-1}\big|_{h=1}=\mathrm{C}^{-1}_f f'(\mathrm{D})^{1/\lambda}\mathrm{C}_{\mathfrak{T}f}\left[ [ (1+\lambda\theta)\omega'(\mathrm{D})^{-1} ]^2 x\frac{ h(1+\beta\theta)}{(1+\lambda \theta)(1+\kappa \theta)}\right] \mathrm{C}_{\mathfrak{T}f}^{-1}f'(\mathrm{D})^{-1/\lambda} \mathrm{C}_f=\\
&=(1+\lambda\theta)^2 ~\mathrm{C}^{-1}_f f'(\mathrm{D})^{1/\lambda}\mathrm{C}_{\mathfrak{T}f}\left[ x\frac{ h(1+\beta\theta)}{(1+\lambda \theta)(1+\kappa \theta)}\right] \mathrm{C}_{\mathfrak{T}f}^{-1}f'(\mathrm{D})^{-1/\lambda} \mathrm{C}_f=\\
&~~~~~~~~~~~~= (1+\lambda\theta)^2 r\left[x \frac{1}{1+\lambda\theta}+a+\tfrac{1}{4}\lambda a^2 \frac{2+\lambda\theta}{1+\lambda(\theta+1)}\mathrm{D}\right] +\\
&~~~~~~~~~~~~~~~~~~+(1+\lambda\theta)^2(1-r)\left[x \frac{1}{1+\kappa\theta} + \frac{1}{2}a \frac{2+\lambda\theta}{1+\kappa\theta} +\frac{1}{2}a \frac{\lambda\theta}{1+\kappa(\theta-1)}+\tfrac{1}{4} \lambda a^2 \frac{2+\lambda \theta}{1+\kappa\theta}\mathrm{D}\right]
\end{align*}
Thus the [$h=1$]-part of the associated dual raising operator is also multiplied and is of the form
$$
\frac{\mathcal{H}_{\theta+c}}{\mathcal{H}_{c}}\frac{(c+\theta)_\theta}{\theta!}\mathrm{U}_{\mathcal{W}^*}(c)\frac{\theta!}{(c+\theta)_\theta}\frac{\mathcal{H}_c}{\mathcal{H}_{\theta+c}}\bigg|_{h=1}=(1+\lambda(\theta+c))^2\frac{\mathcal{F}_{\theta+c}(\beta)}{\mathcal{F}_{c}(\beta)}\frac{(c+\theta)_\theta}{\theta!}\mathrm{U}_{\mathcal{J}^*}(c,\beta)\frac{\theta!}{(c+\theta)_\theta}\frac{\mathcal{F}_{c}(\beta)}{\mathcal{F}_{\theta+c}(\beta)}
$$
where $\mathrm{U}_{\mathcal{J}^*}(c)$ is the obtained dual raising operator of associated Jacobi family with parameter $\beta$.
\newpage
\noindent (so that multiplication by $a_\theta$ for associated family simply turns into multiplication by $a_{\theta+c}$ after conjugation). At the same time the [$h=0$]-part is simply
$$
\frac{\mathcal{H}_{\theta+c}}{\mathcal{H}_{c}}\frac{(c+\theta)_\theta}{\theta!}\mathrm{U}_{\mathcal{W}^*}(c)\frac{\theta!}{(c+\theta)_\theta}\frac{\mathcal{H}_c}{\mathcal{H}_{\theta+c}}\bigg|_{h=0}=\frac{\mathcal{F}_{\theta+c}(\widetilde{\beta})}{\mathcal{F}_{c}(\widetilde{\beta})}\frac{(c+\theta)_\theta}{\theta!}\mathrm{U}_{\mathcal{J}^*}(c,\widetilde{\beta})\frac{\theta!}{(c+\theta)_\theta}\frac{\mathcal{F}_{c}(\widetilde{\beta})}{\mathcal{F}_{\theta+c}(\widetilde{\beta})}
$$
Now recall that we know the result of conjugation with $\mathrm{C}_{\mathfrak{T}f}^{-1} f'(\mathrm{D})^{-c-1/\lambda}\mathrm{C}_f$ for associated Jacobi family with any parameter $\beta$. Thus to establish the whole answer we need to know the result of the conjugation of the operator $(1+\lambda(\theta+c))^2$. But this is almost the same as in case $c=0$. For any $f$ with $\omega=(\mathfrak{T}f)^{inv}$ (and in the last line we substitute explicit $f$: $f'(y)=(1+\tfrac{1}{2}\lambda af(y))^{2}$, $\varphi =f^{inv}=y(1+\tfrac{1}{2}\lambda a y)^{-1}$, $f/f'=y-\tfrac{1}{2}\lambda ay^2$, $\omega = (f/f')^{inv} = \tfrac{1-\sqrt{1-2\lambda a y}}{\lambda a}$, i.e. $\omega'(y)^{-2}=1-2\lambda a y$).
\begin{align*}
\mathrm{C}_{\mathfrak{T}f}^{-1}f'(\mathrm{D})^{-c-\tfrac{1}{\lambda}}\mathrm{C}_f \left[ (1+\lambda(\theta+c))^2 \right]& \mathrm{C}^{-1}_f f'(\mathrm{D})^{c+\tfrac{1}{\lambda}}\mathrm{C}_{\mathfrak{T}f}=\\
=&\mathrm{C}_{\mathfrak{T}f}^{-1} f'(\mathrm{D})^{-c-\tfrac{1}{\lambda}} \left[ (1+\lambda(c+x\mathfrak{T}f(\mathrm{D}))^2 \right] f'(\mathrm{D})^{c+\tfrac{1}{\lambda}}\mathrm{C}_{\mathfrak{T}f}=\\
=&\mathrm{C}_{\mathfrak{T}f}^{-1} \left[ ((1+\lambda c)(\mathfrak{T}f)'(\mathrm{D})+\lambda x\mathfrak{T}f(\mathrm{D}))^2 \right]\mathrm{C}_{\mathfrak{T}f}=\\
=&(1+\lambda(\theta+c))\frac{1}{\omega'(\mathrm{D})}(1+\lambda(\theta+c))\frac{1}{\omega'(\mathrm{D})}=\\
=&(1+\lambda(\theta+c))^2-\tfrac{2a\lambda^2}{\kappa}(1+\lambda(\theta+c))(1+\kappa(\theta+c))\mathrm{D}
\end{align*}
We may thus write the following expansion
\begin{footnotesize}
\begin{align*}
&\mathrm{C}_{\mathfrak{T}f}^{-1}f'(\mathrm{D})^{-c-\tfrac{1}{\lambda}}\mathrm{C}_f \left[\frac{\mathcal{H}_{\theta+c}}{\mathcal{H}_{c}}\frac{(c+\theta)_\theta}{\theta!}\mathrm{U}_{\mathcal{W}^*}(c)\frac{\theta!}{(c+\theta)_\theta}\frac{\mathcal{H}_c}{\mathcal{H}_{\theta+c}} \right] \bigg|_{h=1}\mathrm{C}^{-1}_f f'(\mathrm{D})^{c+\tfrac{1}{\lambda}}\mathrm{C}_{\mathfrak{T}f}=\\
&= \left[ (1+\lambda(\theta+c))\frac{1}{\omega'(\mathrm{D})} \right]^2 \frac{(c+\theta)(1+\beta(c+\theta-1))}{(1+\lambda(c+\theta-1))(1+\kappa(c+\theta-1))}f'(\omega(\mathrm{D}))^{1-c-\tfrac{1}{\lambda}} x\frac{1}{1+\theta} f'(\omega(\mathrm{D}))^{c+\tfrac{1}{\lambda}-1}\\
&=\frac{(1+\lambda(\theta+c))^2(c+\theta)(1+\beta(c+\theta-1))}{(1+\lambda(c+\theta-1))(1+\kappa(c+\theta-1))}f'(\omega(\mathrm{D}))^{1-c-\tfrac{1}{\lambda}} x\frac{1}{1+\theta} f'(\omega(\mathrm{D}))^{c+\tfrac{1}{\lambda}-1}-\frac{2a\lambda^2}{\kappa}(1+c+\theta)(1+\beta(c+\theta))
\end{align*}
\end{footnotesize}
And for $h=0$ we have
\begin{align*}
\mathrm{C}_{\mathfrak{T}f}^{-1}f'(\mathrm{D})^{-c-\tfrac{1}{\lambda}}&\mathrm{C}_f \left[\frac{\mathcal{H}_{\theta+c}}{\mathcal{H}_{c}}\frac{(c+\theta)_\theta}{\theta!}\mathrm{U}_{\mathcal{W}^*}(c)\frac{\theta!}{(c+\theta)_\theta}\frac{\mathcal{H}_c}{\mathcal{H}_{\theta+c}} \right] \bigg|_{h=0}\mathrm{C}^{-1}_f f'(\mathrm{D})^{c+\tfrac{1}{\lambda}}\mathrm{C}_{\mathfrak{T}f}=\\
&=\frac{(c+\theta)(1+\widetilde{\beta}(c+\theta-1))}{(1+\lambda(c+\theta-1))(1+\kappa(c+\theta-1))}f'(\omega(\mathrm{D}))^{1-c-\tfrac{1}{\lambda}} x\frac{1}{1+\theta} f'(\omega(\mathrm{D}))^{c+\tfrac{1}{\lambda}-1}
\end{align*}
So that we may finally write down the full expansion
\begin{align*}
\mathrm{C}_{\mathfrak{T}f}^{-1}&f'(\mathrm{D})^{-c-\tfrac{1}{\lambda}}\mathrm{C}_f \left[\frac{\mathcal{H}_{\theta+c}}{\mathcal{H}_{c}}\frac{(c+\theta)_\theta}{\theta!}\mathrm{U}_{\mathcal{W}^*}(c)\frac{\theta!}{(c+\theta)_\theta}\frac{\mathcal{H}_c}{\mathcal{H}_{\theta+c}} \right]\mathrm{C}^{-1}_f f'(\mathrm{D})^{c+\tfrac{1}{\lambda}}\mathrm{C}_{\mathfrak{T}f}=\\
&=(c+\theta)\frac{\mathcal{H}_{\theta+c}}{\mathcal{H}_{\theta+c-1}}f'(\omega(\mathrm{D}))^{1-c-\tfrac{1}{\lambda}} x\frac{1}{1+\theta} f'(\omega(\mathrm{D}))^{c+\tfrac{1}{\lambda}-1} - \genfrac{}{}{0.5pt}{1}{2ah\lambda^2}{\kappa}(1+c+\theta)(1+\beta(c+\theta))
\end{align*}
Now we again denote $k(y)=f'(\omega(y))^{c-1+1/\lambda}$ and $H_\theta = (-1)^\theta(-c)_\theta=(c+\theta-1)_\theta$, so that $H_{\theta+1}/H_\theta = c+\theta$. We also define $\tilde{\mathrm{C}}_2(c)$: $\tilde{\mathrm{C}}_2^{-1}(c)x \tilde{\mathrm{C}}_2(c) = x - 2ah\lambda^2\kappa^{-1}(1+(\theta+c))(1+\beta(\theta+c))$ (so that $\tilde{\mathrm{C}}_2(0)=\tilde{\mathrm{C}}_2$). Now let's look at the result closely. We have obtained an operator of the form
$$
(c+\theta)\frac{\mathcal{H}_{\theta+c}}{\mathcal{H}_{\theta+c-1}}k(\mathrm{D})^{-1} x\frac{1}{1+\theta}k(\mathrm{D}) -  2ah\lambda^2\kappa^{-1}(1+(\theta+c))(1+\beta(\theta+c))
$$
By long division lemma in terms of $\delta=1-x(1+\theta)^{-1}\mathrm{D}$ this operator is equal to
\begin{align*}
(c+\theta)\frac{\mathcal{H}_{\theta+c}}{\mathcal{H}_{\theta+c-1}}&\left(x\frac{1}{1+\theta}-\delta (\mathrm{L}k)(\mathrm{D})\right)-  2ah\lambda^2\kappa^{-1}(1+(\theta+c))(1+\beta(\theta+c))=\\
&=\frac{\mathcal{H}_{\theta+c}}{\mathcal{H}_{c}}\frac{(c+\theta)_\theta}{\theta!}\tilde{\mathrm{C}}_2^{-1}(c)x \tilde{\mathrm{C}}_2(c)\frac{\theta!}{(c+\theta)_\theta}\frac{\mathcal{H}_c}{\mathcal{H}_{\theta+c}} - c\frac{\mathcal{H}_{c}}{\mathcal{H}_{c-1}}\delta(\mathrm{L}k)(\mathrm{D})
\end{align*}
(recall that $\delta$ annihilates operators, so $m_\theta \delta = \delta m_\theta = m_0 \delta$).
\newpage
\noindent Then conjugation with $\theta!^{-1}\tilde{\mathrm{C}}_2(c)\theta!(c+\theta)_\theta^{-1}\mathcal{H}_c\mathcal{H}_{\theta+c}^{-1}$ leads to the following operator
\begin{align*}
x\frac{1}{1+\theta}-c\frac{\mathcal{H}_{c}}{\mathcal{H}_{c-1}}&\delta\left(\theta!\overline{\tilde{\mathrm{C}}_2(c)^{-1}}\theta!^{-1}(c+\theta)_\theta\mathcal{H}_c^{-1}\mathcal{H}_{\theta+c}\mathrm{L}\cdot k\right)(\mathrm{D})=\\
&x\frac{1}{1+\theta}-\delta\left(\theta!\overline{\tilde{\mathrm{C}}_2(c)^{-1}}\theta!^{-1}\mathrm{L}(c+\theta-1)_\theta\mathcal{H}_{c-1}^{-1}\mathcal{H}_{\theta+c-1}\cdot k\right)(\mathrm{D})
\end{align*}
We now denote the series $s(c,y)\coloneqq 1+y\theta!\overline{\tilde{\mathrm{C}}_2(c)^{-1}}\theta!^{-1}\mathrm{L}(c+\theta-1)_\theta\mathcal{H}_{c-1}^{-1}\mathcal{H}_{\theta+c-1}\cdot f'(\omega(y))^{c-1+1/\lambda}$. So we obtain the following result:
\begin{align*}
\mathrm{C}_{\mathfrak{T}f}^{-1}f'(\mathrm{D})^{-c-\tfrac{1}{\lambda}}\mathrm{C}_f &\left[\frac{\mathcal{H}_{\theta+c}}{\mathcal{H}_{c}}\frac{(c+\theta)_\theta}{\theta!}\mathrm{U}_{\mathcal{W}^*}(c)\frac{\theta!}{(c+\theta)_\theta}\frac{\mathcal{H}_c}{\mathcal{H}_{\theta+c}} \right]\mathrm{C}^{-1}_f f'(\mathrm{D})^{c+\tfrac{1}{\lambda}}\mathrm{C}_{\mathfrak{T}f}=\\
&=\frac{\mathcal{H}_{\theta+c}}{\mathcal{H}_{c}}\frac{(c+\theta)_\theta}{\theta!}\tilde{\mathrm{C}}_2^{-1}(c)\theta!s(c,\mathrm{D})^{-1}x\frac{1}{1+\theta}s(c,\mathrm{D}) \theta!^{-1}\tilde{\mathrm{C}}_2(c)\frac{\theta!}{(c+\theta)_\theta}\frac{\mathcal{H}_c}{\mathcal{H}_{\theta+c}}
\end{align*}
Hence we have the following explicit representation for the $\mathrm{G_P}$-operator of associated Wilson family
\begin{align*}
\mathrm{G}_{\mathcal{W}}(c)=\theta! s(c,\mathrm{D}) \theta!^{-1}\tilde{\mathrm{C}}_2(c)\frac{\theta!}{(c+\theta)_\theta}\frac{\mathcal{H}_c}{\mathcal{H}_{\theta+c}}\mathrm{C}_{\mathfrak{T}f}^{-1}f'(\mathrm{D})^{-c-\tfrac{1}{\lambda}}\mathrm{C}_f  \frac{\mathcal{H}_{\theta+c}}{\mathcal{H}_{c}}\frac{(c+\theta)_\theta}{\theta!}
\end{align*}
Notice that we did not symmetrize our notations, thus we may not correctly rewrite the moment generating function as a fraction of two hypergeometric series of higher order. The reason is that we do not know explicit factorization of polynomial $(1+\lambda \theta)(1+\kappa \theta)\mathcal{H}_{\theta+1}/\mathcal{H}_\theta=(h(1+\beta\theta)(1+\lambda(\theta+1))^2+(1-h)(1+\widetilde{\beta}\theta))$, which is a polynomial in $\theta$ of degree 3, so that it requires some work. 'Hypergeometricity' of the result follows from the fact that if $\tilde{\mathrm{C}}_\ell^{-1}x \tilde{\mathrm{C}}_\ell = x - \ell_\theta $, then in addition to the factorization of corresponding polynomials $\tilde{\mathrm{C}}_\ell \cdot x^n = \prod_{k=0}^{n-1} (x+\ell_k)$ we also have factorization
$$
\theta! \overline{\tilde{\mathrm{C}}_\ell^{-1}}\theta!^{-1} \cdot y^n =\frac{y^n}{(1+y\ell_0)~...~(1+y\ell_n)}
$$
since
$$
\theta! \overline{\tilde{\mathrm{C}}_\ell^{-1}}\theta!^{-1} (\mathrm{L}-\ell_\theta)=\mathrm{L}~\theta! \overline{\tilde{\mathrm{C}}_\ell^{-1}}\theta!^{-1} 
$$
\noindent\rule{\textwidth}{1.2pt}
~\\
\noindent [Second comment on established polynomials]. In summary, what we did is that we expanded the Jacobi-Wilson generators in the following manner
\begin{align*}
&\frac{1+u_1\theta}{1+\lambda(\theta-1)} f'(\omega(\mathrm{D}))^{1-1/\lambda} x \frac{1}{1+\theta} f'(\omega(\mathrm{D}))^{-1+1/\lambda}\\
&\frac{1+u_2\theta}{1+\kappa(\theta-1)} f'(\omega(\mathrm{D}))^{1-1/\lambda} x \frac{1}{1+\theta} f'(\omega(\mathrm{D}))^{-1+1/\lambda}\\
\frac{(1+\lambda\theta)^2(1+u_3\theta)}{1+\lambda(\theta-1)} &f'(\omega(\mathrm{D}))^{1-1/\lambda} x \frac{1}{1+\theta} f'(\omega(\mathrm{D}))^{-1+1/\lambda}-\tfrac{2\lambda^2a}{\kappa}(1+\kappa\theta)(1+u_3(1+\theta))\\
\frac{(1+\lambda\theta)^2(1+u_4\theta)}{1+\kappa(\theta-1)} &f'(\omega(\mathrm{D}))^{1-1/\lambda} x \frac{1}{1+\theta} f'(\omega(\mathrm{D}))^{-1+1/\lambda}-\tfrac{2\lambda^2a}{\kappa}(1+\lambda\theta)(1+u_4(1+\theta))
\end{align*}
And noticed that the diagonalization of an arbitrary linear combination of them is well-understood. Unfortunately, much more interesting cases, such as
\begin{align*}
f'(\omega&(\mathrm{D}))^{\tfrac{1}{\sigma}-\tfrac{1}{\lambda}}x\frac{1}{1+\sigma\theta}f'(\omega(\mathrm{D}))^{\tfrac{1}{\lambda}-\tfrac{1}{\sigma}} ~~~~~~~~~~~~~~~~~ x\frac{1}{\omega'(\mathrm{D})}\frac{1}{1+\lambda\theta}\omega'(\mathrm{D})\frac{1}{1+\lambda\theta}\\
&x\frac{1}{1+\lambda\theta}\omega'(\mathrm{D})\frac{1}{1+\lambda\theta}\omega'(\mathrm{D})\frac{1}{1+\lambda\theta}~~~~~~~~~~~x\frac{1}{1+\lambda\theta}+t(1+\lambda\theta)\frac{1}{\omega'(\mathrm{D})}
\end{align*}
(so that with further division of the coefficients in three-term recurrence relations or deformation of the powers) or their linear combinations seem harder to deal with, but as raising operators all of them obviously make perfect sense, so that the corresponding orthogonal polynomials may still be calculated by hand.\\
\noindent\rule{\textwidth}{1.2pt}


\newpage
 \begin{center}
\textbf{Conclusion}
\end{center}
Such an approach, introduced in this paper, allows us to naturally vary the base power series $f$. Consider, for example the Sheffer $[n+1]$-term recurrence relation, corresponding to polynomial family with the dual raising operator of the form
$$
\mathrm{U_{P^*}}=(1+\tfrac{1}{n}\lambda a\mathrm{D})^{n/\lambda} x (1+\tfrac{1}{n}\lambda a\mathrm{D})^{n-n/\lambda}
$$
where we still have $\mathrm{G_P}=f'(\mathrm{D})^{-1/\lambda} \mathrm{C}_f$ but now $f'(y)=(1+\tfrac{1}{n}\lambda af(y))^{n}$, so that the derivative of $\varphi = f^{inv}$ is defined as $\varphi'(y)=(1+\tfrac{1}{n}\lambda ay)^{-n}$. Now fix $\lambda_0 = \lambda$, and define $\lambda_k$ by the relation $\tfrac{1}{\lambda_k}-\tfrac{1}{\lambda}=\tfrac{k}{n}$ for $k=1,...,n-1$. Consider now the following set of dual raising operators
\begin{align*}
\varphi'(\mathrm{D})^{-\tfrac{1}{\lambda}} ~x \varphi'(\mathrm{D})^{\tfrac{k}{n}-1+\tfrac{1}{\lambda}}\frac{1}{1+\lambda_k \theta} \varphi'(\mathrm{D})^{-\tfrac{k}{n}} &=  (1+\tfrac{1}{n}\lambda a\mathrm{D})^{\tfrac{n}{\lambda}} x (1+\tfrac{1}{n}\lambda a\mathrm{D})^{n-k-\tfrac{n}{\lambda}} \frac{1}{1+\lambda_k \theta} (1+\tfrac{1}{n}\lambda a\mathrm{D})^k
\end{align*}
Then for any $k=0,1,...,n-1$ these operators are $[n+1]$-term. Now exactly as in case of Jacobi orthogonal family we notice that the operator on the left may be rewritten as
$$
\mathrm{C}^{-1}_f f'(\mathrm{D})^{1/\lambda} \mathrm{C}_{\mathfrak{T}f} ~x \frac{1}{1+\lambda_k \theta}~\mathrm{C}_{\mathfrak{T}f}^{-1} f'(\mathrm{D})^{-1/\lambda} \mathrm{C}_f=\varphi'(\mathrm{D})^{-\tfrac{1}{\lambda}} ~x \varphi'(\mathrm{D})^{-1+\tfrac{1}{\lambda_k}}\frac{1}{1+\lambda_k \theta} \varphi'(\mathrm{D})^{\tfrac{1}{\lambda}-\tfrac{1}{\lambda_k}}
$$
Any such operator is $[n+1]$-term, hence their linear combination is $[n+1]$-term, so that for any $t_i$: $\sum_0^{n-1} t_i =1$ the following operator is $[n+1]$-term.
$$
\mathrm{C}^{-1}_f f'(\mathrm{D})^{1/\lambda} \mathrm{C}_{\mathfrak{T}f} \left[\sum_{k=0}^{n-1} x \frac{t_k}{1+\lambda_k \theta}\right]\mathrm{C}_{\mathfrak{T}f}^{-1} f'(\mathrm{D})^{-1/\lambda} \mathrm{C}_f
$$
Hence we may define an invertible operator $\mathcal{F}_\theta$, such that $\mathcal{F}_{\theta+1}/\mathcal{F}_\theta=\sum_{k=0}^{n-1} t_k(1+\lambda_k \theta)^{-1}$, $\mathcal{F}_0=1$. Then the direct analogue of Jacobi polynomials in this case is the polynomial family $\mathrm{Q}$ defined as
$$
\mathrm{G_Q}=\mathcal{F}_\theta^{-1} \mathrm{C}_{\mathfrak{T}f}^{-1} f'(\mathrm{D})^{-1/\lambda} \mathrm{C}_f \mathcal{F}_\theta
$$
and the calculations above show that these polynomials for any parameters $t_i$ satisfy $[n+1]$-term recurrence relation. Now the natural question arises, if there is a way to obtain $2n$ $[n+1]$-generators as in case of Wilson polynomials. The author doesn't know, if there is a way to do that, since the series $\omega=(f/f')^{inv}$ in this case is more complicated. But there is at least one additional generator, corresponding to the case of Hahn polynomials. Notice that
\begin{align*}
\varphi = &\frac{(1+\tfrac{1}{n}\lambda ay)^{1-n}-1}{\frac{1-n}{n}\lambda a} ~~~\Rightarrow~~~ f=\frac{(1+\tfrac{1-n}{n}\lambda ay)^{\tfrac{1}{1-n}}-1}{\frac{1}{n}\lambda a}~~~\Rightarrow~~~\mathfrak{T}f=\frac{(1+\tfrac{1-n}{n}\lambda ay)^{\tfrac{1}{1-n}}-1}{\frac{1}{n}\lambda a(1+\tfrac{1-n}{n}\lambda ay)^{\tfrac{n}{1-n}}}\\
&~~~\Rightarrow~~~ (\mathfrak{T}f)' = 1-n + n(1+\tfrac{1-n}{n}\lambda ay)^{\tfrac{1}{n-1}} ~~~\Rightarrow~~~ (\mathfrak{T}f)'(\varphi)=\frac{1+\frac{1-n}{n}\lambda a y}{1+\frac{1}{n}\lambda a y}
\end{align*}
We thus have the following algebraic relation for $\mathfrak{T}f$ and its inverse $\omega(y)$
$$
(1-(\mathfrak{T}f)')((\mathfrak{T}f)'+n-1)^{n-1}=n^{n-1}\lambda a \mathfrak{T}f~~~~\Rightarrow~~~~\left(1-\omega'(y)^{-1}\right)\left(\omega'(y)^{-1}+n-1\right)^{n-1} = n^{n-1}\lambda ay
$$
Now write $P(y)=(n-1)^{n-1}+(y-1)(y+n-1)^{n-1}$. This is a polynomial of degree $n$ and $P(0)=0$,\linebreak hence $P(y)/y$ is a polynomial of degree $n-1$. Also $P(\omega'(y)^{-1})= (n-1)^{n-1}-n^{n-1}\lambda ay$. We now claim that the following operator is also $[n+1]$-term
\[
\mathrm{C}^{-1}_f f'(\mathrm{D})^{1/\lambda} \mathrm{C}_{\mathfrak{T}f} ~\left[x (n-1)^{n-1}-n^{n-1}\lambda a\theta\right]~\mathrm{C}_{\mathfrak{T}f}^{-1} f'(\mathrm{D})^{-1/\lambda} \mathrm{C}_f
\]
And that is indeed true since it may be rewritten as
$$
\mathrm{C}^{-1}_f f'(\mathrm{D})^{1/\lambda} \mathrm{C}_{\mathfrak{T}f} \left[x P(\omega'(\mathrm{D})^{-1})\right]\mathrm{C}_{\mathfrak{T}f}^{-1} f'(\mathrm{D})^{-1/\lambda} \mathrm{C}_f = \varphi'(\mathrm{D})^{-\tfrac{1}{\lambda}} x \frac{P((\mathfrak{T}f)'(\varphi(\mathrm{D})))}{(\mathfrak{T}f)'(\varphi(\mathrm{D}))} \varphi'(\mathrm{D})^{\tfrac{1}{\lambda}-1}
$$
Now, for any $m=0,...,n-1$ the operator
$$
\varphi'(\mathrm{D})^{-\tfrac{1}{\lambda}} x (\mathfrak{T}f)'(\varphi(\mathrm{D}))^m \varphi'(\mathrm{D})^{\tfrac{1}{\lambda}-1} = (1+\tfrac{1}{n}\lambda a\mathrm{D})^{n/\lambda} x (1+\tfrac{1-n}{n}\lambda a\mathrm{D})^{m}(1+\tfrac{1}{n}\lambda a\mathrm{D})^{n-m-n/\lambda}
$$
is clearly $[n+1]$-term. Since $P(y)/y$ is a polynomial of degree $n-1$, the result follows. Thus we have the following extension in this case. For any $t_i$: $\sum_0^n t_i =1$ the following operator is $[n+1]$-term after conjugation with $\mathrm{C}^{-1}_f f'(\mathrm{D})^{1/\lambda} \mathrm{C}_{\mathfrak{T}f}$:
$$
\left[\sum_{k=0}^{n-1} x \frac{t_k}{1+\lambda_k \theta}\right] + t_n \left[\vphantom{\sum_0^n}x - \frac{n^{n-1}}{(n-1)^{n-1}}\lambda a\theta\right]
$$
And this operator is well-understood. Define $\mathcal{Q}_\theta$: $\mathcal{Q}_{\theta+1}/\mathcal{Q}_\theta = t_n + \sum_{k=0}^{n-1} t_k(1+\lambda_k \theta)^{-1}$, $\mathcal{Q}_0=1$, the constant $\mathcal{C}=-t_n\lambda an^{n-1}(n-1)^{1-n}$ and consider the function $\Delta_{\mathcal{C}}=(e^{\mathcal{C} y}-1)/\mathcal{C}$. Then the operator above may be rewritten as
$$
\left[\sum_{k=0}^{n-1} x \frac{t_k}{1+\lambda_k \theta}\right] + t_n \left[\vphantom{\sum_0^n}x - \frac{n^{n-1}}{(n-1)^{n-1}}\lambda a\theta\right] = \mathcal{Q}_\theta \mathrm{C}_{\Delta_{\mathcal{C}}}^{-1} x \mathrm{C}_{\Delta_{\mathcal{C}}}^{\vphantom{-1}} \mathcal{Q}_\theta^{-1}
$$
Thus the new family of monic polynomials $\mathrm{Q}$ defined as
$$
\mathrm{G_Q}=\mathrm{C}_{\Delta_{\mathcal{C}}} \mathcal{Q}_\theta^{-1} \mathrm{C}_{\mathfrak{T}f}^{-1} f'(\mathrm{D})^{-1/\lambda} \mathrm{C}_f \mathcal{Q}_\theta ~~~~~~~~~~~~~ p_n(x)=\mathrm{G_Q}\cdot x^n
$$
satisfies $[n+1]$-term recurrence relation. It is an interesting question, if this family may be extended further, does it have any important properties or any useful applications.
~\\
~\\

\newpage
\begin{center}
\textsc{Appendix A.} \textit{Why $p_{-1}(y)\neq 0$}
\end{center}

\noindent [Possible duality]. Consider the general case of recurrence relations of an arbitrary order. Write
$\mathrm{U}_{\mathrm{P}^*}=x+\sum_{k\geqslant 0} \frac{1}{k!} \mathrm{D}^k a_\theta^k ~\iff~ x=\mathrm{U}_{\mathrm{P}}+\sum_{k\geqslant 0}  \frac{1}{k!} \mathrm{D}_\mathrm{P}^k a_{\theta_{\mathrm{P}}}^k$. In other words we have the following relation
$$
xp_n(x)=p_{n+1}(x)+\sum_{k=0}^\infty \binom{n}{k} a^k_n p_{n-k}(x)
$$
Writing the latter in terms of derivative operator for $f_\ell(y)$:
$$
\mathfrak{d}=\mathfrak{d}_\mathrm{P} + a_{\mathfrak{u}_\mathrm{P}\mathfrak{d}_{\mathrm{P}}}^0 + a_{\mathfrak{u}_\mathrm{P}\mathfrak{d}_{\mathrm{P}}}^1 \mathfrak{u}_\mathrm{P} + a_{\mathfrak{u}_\mathrm{P}\mathfrak{d}_{\mathrm{P}}}^2 \mathfrak{u}_\mathrm{P}^2 \frac{1}{2!}+a_{\mathfrak{u}_\mathrm{P}\mathfrak{d}_{\mathrm{P}}}^3 \mathfrak{u}_\mathrm{P}^3 \frac{1}{3!}+...
$$
We obtain the following recurrence relation for $\ell!F_\ell(x)=\theta!\cdot f_\ell(x)$
$$
\mathrm{L}\cdot F_{n}(y)= \delta_n^0 F_{n-1}(y) + \sum_{k=0}^\infty \binom{n+k}{k} a_{n+k}^k F_{n+k}(y)
$$
And here something interesting actually happens. This formalism allows us to define $F_{-1}(y)$, and it happens in a very natural way. Setting $F_{-1}(y)=y^{-1}$ we obtain that for all $n$ by definition of $\mathrm{L}$:
$$
y^{-1}F_{n}(y) = F_{n-1}(y) + \sum_{k=0}^\infty \binom{n+k}{k} a_{n+k}^k F_{n+k}(y)
$$
(Notice that this approach does not allow to define $f_{-1}(y)$ by contrast, since $\theta!$ is not defined correctly in negative domain). So that we naturally define
$$
\Hat{p}_{-n-1}(y) \coloneqq y^{-1}F_{n}(y^{-1}) ~~~\Rightarrow~~~ n!\Hat{p}_{-n-1}(y) = \int_{0}^\infty f_{n}(t) e^{-yt} dt
$$
Moreover, if the initial $a_\theta$ were polynomials in $\theta$, we can now naturally define
$$
\boxed{
\Hat{a}_\theta^k = (-1)^k a_{k-1-\theta}^k}
$$
in order to define the dual \textit{polynomials} $\Hat{p}_n(x)$ satisfying recurrence relation
$$
x \Hat{p}_n(x) = \Hat{p}_{n+1}(x) + \sum_{k=0}^\infty \binom{n}{k} \Hat{a}^k_n \Hat{p}_{n-k}(x)
$$
And this recurrence relation now holds for all integer indices, including the negative ones.\\
~\\
\textit{Remark}: Here I have to explain why I think that this topic is, well, not quite important, but at least should be treated carefully. In the most part of the textbooks on orthogonal polynomials and related topics the authors set by default $p_{-1}(x)=0$ in three-term recurrence relations. And that contradicts what I think should be done. Since we work with monic polynomials, $p_{-1}$ should be naturally defined as a series of the form $x^{-1}+x^{-2}\mathbb{C}[[x^{-1}]]$. We see that it cannot be correctly done in the general case for the obvious reason, but the case of polynomial coefficients in recurrence relations already gives us some insights towards this direction. Actually exactly this procedure allows us to define the negative-index 'polynomials' of binomial type (just take this definition of the dual family and figure out the general formula for the dual of the Sheffer family). It also allows us at least to feel, how polynomials behave under the variation of $f_n(y) \rightarrow f_n(y) e^{-any}$. We just take the dual, then consider the new family $\Hat{p}_n(x-a(n+1))$ and then take the dual again. Moreover this type of duality sends orthogonal polynomials to orthogonal ones, and then if we define the dual moments, the following identity holds
$$
p_{-1}(x)=\sum_{n=0}^\infty \frac{n!\Hat{\mathcal{B}}_n}{\Hat{p}_n(x)\Hat{p}_{n+1}(x)}
$$

\end{document}